\documentclass[12pt]{article}
\usepackage{latexsym,amsmath,amssymb,amsbsy,amstext,amscd,amsfonts}
\usepackage{color}
\usepackage{pictex,epsfig}
\usepackage[notcite,notref]{}

\newtheorem{theorem}{\noindent {\rm \bf Theorem}}[section]

\begin{document}

\begin{center}
{\bf Exact solution to a nonlinear heat conduction problem in doubly
periodic 2D composite materials}
 \end{center}

\begin{center}
{ D. Kapanadze$^{\sharp,}$, G. Mishuris$^{*}$, E.
Pesetskaya$^{\sharp,}$\footnote{Corresponding author.}}
 \end{center}

\begin{center}
$^{*}$ Department of Mathematics and Physics, Aberystwyth
University, UK \\
$^{\sharp}$ A. Razmadze Mathematical Institute, Tbilisi State University, Georgia
 \end{center}

 {\bf Abstract }

An analytic solution to a stationary heat conduction problem in 2D unbounded doubly periodic composite
materials with temperature dependent conductivities
of their components is given. Corresponding nonlinear boundary value
problem is reduced a Laplace equation with nonlinear transmission conditions.
For special relationships between the conductivity
coefficients of the matrix and inclusions, the problem is transformed
to fully linear boundary value problem for doubly periodic analytic functions.
This allows to reconstruct the solution of the originally nonlinear composite and
to find its effective properties. The results are illustrated by numerical examples.

\section{Introduction }

In the paper, a stationary nonlinear heat conduction problem in a 2D unbounded doubly periodic composite with inclusions is considered.
Conductive properties of the inclusions assumed to be proportional to that of the matrix. For such nonlinear composite, an exact solution is constructed and the effective conductivity of the composite is discussed in details.

The theory and technique for solution of the linear boundary value problems for 2D doubly periodic
composite materials with constant conductivities of their components are well developed
in discussed in details in \cite{BeMi01}, \cite{KPM}, \cite{Kusch} etc.

In case of nonlinear heat conduction problem, two classes of the problems can be identified. The first one is when the
material parameters depend on gradient of the temperature, and the second one is when the parameters are functions of the temperature itself.
It was shown in \cite{MishSev}, where the homogenization procedure for a random composite with conductivities dependent on temperature have been developed basing on the classical approaches,
that the latter problem  is more difficult in comparison with the former.
Among others, the authors also proved that the Eshelby inclusion approach is not valid when the material parameters are functions of temperature.

The most complete analysis of the nonlinear heat problem in periodic composite
was done in \cite{GalTelTok}.
The authors used rigorous asymptotic homogenization technique for periodic micro-heterogeneous (see for example \cite{BakPan,BenLions}) and
shown that the leading asymptotic term can be constructed from solution of a homogenization problem for a specific linear heterogeneous composite.
To evaluate the average properties of the composites, a local (RVE level) problem was formulated as an abstract minimization problem.
Hashin-Shtrikman estimates for the material parameters were also given.

Average properties of such composites were evaluated by means of Pad\'{e} approximation
approach in \cite{Tok}, \cite{Tok2}. Some results concerning on effective properties of special 2D doubly periodic porous media was discussed
in \cite{Mit98} and \cite{Pes} on the base of the analytic functions approach.

However, the problem for nonlinear composite materials  where the conductivity dependents on the temperature is far from completeness.
In present paper the first exact solution for the double periodic nonlinear composite is constructed under specific assumptions of the material
properties of the composite components. Namely, we consider the static thermal conductivity problem of
unbounded 2D composite materials with circular disjoint inclusions
geometrically formed doubly periodic structure. We suppose that
each component of the composite is filled in by materials of
different conductivity depending on the temperature. The key point in the analysis is the assumption that
ratios of the component conductivities are independent of the temperature.
A steady (external) flux of a given average intensity flows in a certain direction
within the composite and is not, in general, parallel to orientations of the periodic
cell. The components are coupled together into unique structure
due to so called ideal contact conditions.

The main goal of this work is to determine the
temperature and the flux distributions in analytic
forms and to derive average properties of such composites.
In contrast to the linear problem (\cite{KPM}), the flux is not a doubly periodic function for a nonlinear composite. However,
under that aforemention assumption on the material properties, we prove such properties for such composite. Even in this case,
the temperature is no longer quasi-periodic function.

The paper is organized as follows. Accurate formulation of the problem is given in Section 2.
In Section 3, we reduce the given nonlinear boundary value problem defined by a nonlinear partial differential equation and
linear transmission conditions, describing continuity of the temperature and heat flux along the interfaces between the matrix and inclusions, to
an equivalent nonlinear boundary value problem for Laplace equation with nonlinear transmission conditions. Then we formulate conditions for which the
new problem can be linearized. Although, as it was mentioned above, a solution of the nonlinear problem in periodic structure
cannot be represented by doubly periodic function, we show that some ideas coming from linear
techniques could be effectively used for the solution construction for some class of the composite materials important for applications.

Numerical calculations are performed and
discussed in Section 5. In this section, we present average properties of the composite, compare the results with the formula from \cite{MishSev} evaluated for other type of the composites and discuss the obtained results.

\section{Statement of the problem}

 Let us first describe the geometry of the composites.
We consider a lattice $L$ which is defined by the two fundamental
translation vectors $1$ and $\imath$ (where $\imath^2$ = -1) in the
complex plane $\mathbb{C}\cong \mathbb{R}^2$ of the complex variable
$z=x+\imath y$. Here, the representative cell is the
square
$$ Q_{(0,0)} := \left\{ z= t_1+\imath t_2 \in \mathbb{C}:
-\frac{1}{2}<t_{p}<\frac{1}{2},\, p=1,2\right\}.
$$
Let $\mathcal{E}:= \bigcup\limits_{m_1, m_2} \{m_1+\imath m_2\}$ be
the set of the lattice points, where $m_1, m_2 \in \mathbb{Z}$. The
cells corresponding to the points of the lattice $\mathcal{E}$ are
denoted by
$$
Q_{(m_1,m_2)}=Q_{(0,0)}+m_1+\imath m_2:= \left\{z\in \mathbb{C}:
z-m_1-\imath m_2 \in Q_{(0,0)}\right\}.
$$


It is considered the situation when mutually disjoint disks
(inclusions) of different radii $D_{k}:=\{z \in \mathbb{C}:
|z-a_{k}|<r_{k}\}$ with boundaries $\partial D_{k}:=\{z \in
\mathbb{C}: |z-a_{k}|=r_{k}\} \, (k=1,2,\dots,N)$  are located
inside the cell $Q_{(0,0)}$ and periodically repeated in all cells
$Q_{(m_1,m_2)}$. We denote by
$$
D_0:=Q_{(0,0)}\setminus \left(\bigcup\limits_{k=1}^{N}\,D_{k}\cup
\partial D_{k}\right)
$$
the connected domain obtained by removing of the inclusions from
the cell $Q_{(0,0)}$.
\begin{figure}[h!]
\begin{center}
\resizebox*{11cm}{!}{\includegraphics{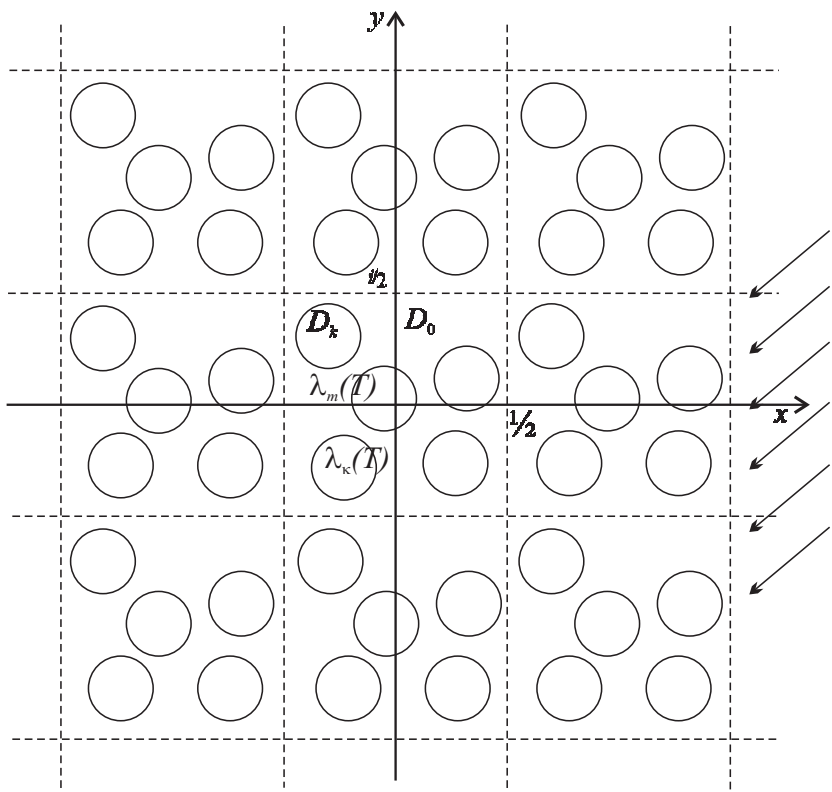}}
\caption{2D double periodic composite with inclusions.}%
\label{fig:pic}
\end{center}
\end{figure}

We investigate the steady state heat conduction problem for
nonlinear composite materials modeling by the above described
geometry, i.e., determination of a distribution of the temperature
$T$ (and/or heat flux $q$) in such composites.

We consider a doubly periodic composite material with matrix
$$
D_{matrix}=\bigcup\limits_{m_1,m_2} \,((D_0 \cup \partial
Q_{(0,0)})+m_1+\imath m_2)
$$
and inclusions
$$
D_{inc}=\bigcup\limits_{m_1,m_2} \bigcup\limits_{k=1}^{N}\,
(D_{k}+m_1+\imath m_2)
$$
occupied by materials of conductivities $\lambda(T)$ and
$\lambda_k(T)$, respectively.

The thermal loading for the composite is described weakly by the flux
given at infinity or more accurately by its intensity $A$. We assume that the flux
is directed $\theta$ which does not coincide, in general, with the
orientation of the periodic cell. According to the
conservation law and the ideal (perfect) contact condition between
the different materials the flux is continuous in the entire
structure. Moreover, as a result of such formulation, the
temperature which is also continuous, as the results of the ideal transmission conditions along the interface between the matrix and inclusions,
possesses non-zero jumps across any cell.

We assume that the conductivities $\lambda(T), \lambda_k(T) \in
\mathcal{C^\infty} ({\mathbb R}), k=1,\dots,N,$ are
continuous positive functions on ${\mathbb R}$ such
that
\begin{equation}
\label{lambda_bound}
 0 <  \lambda(T) < + \infty,
\end{equation}
\begin{equation}
\label{lambda_k_bound} 0 <  \lambda_{k}(T) < + \infty,\;\;\; k =1, \ldots, N.
\end{equation}

We search for steady-state distribution of the temperature
and heat flux within such composite. The problem is equivalent to determination of
the function $T = T(x,y)$ satisfying the
nonlinear differential equation
\begin{equation}\label{1a}
\nabla(\lambda(T)\nabla T)=0,\; (x,y) \in D_{matrix},
\end{equation}
\begin{equation}\label{1b}
\nabla(\lambda_k(T)\nabla T)=0, \; (x,y) \in
\bigcup\limits_{m_1,m_2 \in
\mathbb{Z}}\, D_k+m_1+\imath m_2.
\end{equation}

We assume that the average flux vector of intensity $A$ is directed at an
angle $\theta$ to axis $Ox$ (see Fig. \ref{fig:pic}).
Transmission conditions at the boundary of each cell can be written in the forms
\begin{equation}\label{bound_1}
\lambda(T)T_y|_{\partial Q^{(top)}_{(m_1,m_2)}}=-A
\sin \theta+q^+_1(x+m_1,m_1,m_2),
\end{equation}
\[
\lambda(T)T_y|_{\partial Q^{(bottom)}_{(m_1,m_2)}}=-A
\sin \theta+q^-_1(x+m_1,m_1,m_2),
\]
where $-1/2<x<1/2$ is a local coordinate system connected with each particular cell.

In the perpendicular direction we have similar relationships:
\begin{equation} \label{bound_2}
\lambda(T)T_x|_{\partial Q^{(left)}_{(m_1,m_2)}}=-A
\cos \theta+q^-_2(y+m_2,m_1,m_2),
\end{equation}
\[
\lambda(T)T_x|_{\partial Q^{(right)}_{(m_1,m_2)}}=-A
\cos \theta+q^+_2(y+m_2,m_1,m_2),\quad -1/2<y<1/2.
\]
Note that in general the solution is not periodic even in terms of flux. Thus the unknown flux at the different side of the unit cell
$q_j^\pm(\cdot,m_1,m_2)$, ($j=1,2$) are not the same. However, since there are no sources and sinks in the composite, the following conditions hold true:
\begin{equation} \label{zero_flux}
\int_{-1/2}^{1/2}q^\pm_j(\xi+m_j,m_1,m_2)d\xi=0.
\end{equation}

As a result of the energy conservation law (see also (\ref{bound_1}), (\ref{bound_2}) and (\ref{zero_flux})),
the integral of the heat flux over the cell boundary is equal to zero
\begin{equation}\label{bound_4aa}
\int\limits_{\partial\,
Q_{(m_1,m_2)}}\hspace{-5mm}\lambda(T)\frac{\partial
T}{\partial n}ds=0.
\end{equation}

Finally, we assume that the ideal contact conditions on the boundaries
between the matrix and the inclusions are hold:
\begin{equation}\label{bound_3}
T(t)= T_k(t), \quad t \in \bigcup\limits_{m_1,m_2 \in
\mathbb{Z}}\,(\partial D_k +m_1+\imath m_2),
\end{equation}
\begin{equation}\label{bound_4}
\lambda(T(t))\frac{\partial T(t)}{\partial
n}=\lambda_k(T_k(t))\frac{\partial T_k(t)}{\partial n}, \quad t \in
\bigcup\limits_{m_1,m_2 \in \mathbb{Z}}\,(\partial D_k +m_1+\imath
m_2).
\end{equation}
Here, the vector $n=(n_1,n_2)$ is the outward unit normal vector
to $\partial D_k$; $\frac{\partial}{\partial n}=n_1
\frac{\partial}{\partial x}+n_2 \frac{\partial}{\partial y};\;$
and $ T(t):=\lim\limits_{D_0\ni z \rightarrow t}\, T(z), \;\;
T_k(t):=\lim\limits_{D_k \ni z \rightarrow t}\, T(z)$.


\section{Reformulation of the problem}

To solve the problem, we use the so-called Kirchhoff
transformation \cite{Kirchhoff} (known also as Baiocchi transformation \cite{BaiCap}). Namely,  let us
introduce continuous increasing functions $f: {\mathbb
R} \rightarrow {\mathbb R}, \, f_k: {\mathbb R} \rightarrow
{\mathbb R},\, k=1,\dots,N,$
\begin{equation}
\label{Baiocchi} f(T)=\int\limits_{0}^{T}\,\lambda(\xi)\, d\xi,
\quad f_k(T)=\int\limits_{0}^{T}\,\lambda_k(\xi)\, d\xi, \quad
k=1,\dots,N,
\end{equation}
and perform the following change of the variables:
\begin{equation}
\label{harmonic}
 u(x,y)=f(T(x,y)), \quad u_k(x,y)=f_k(T_k(x,y)), \quad
k=1,\dots,N.
\end{equation}
Note that $f(T)$ and $f_k(T)$ are monotonic functions of temperature and therefore there exist their inverses $f^{-1}$ and $f_k^{-1}$.

By using representations (\ref{Baiocchi}), the equations (\ref{1a}),
(\ref{1b}) are transformed to the Laplace equations (see, e.g.,
\cite{BaiCap})
\begin{equation}\label{Lapl_a}
\Delta u=0,\quad (x,y) \in D_{matrix},
\end{equation}
\begin{equation}\label{Lapl_b}
\Delta u_k=0, \quad (x,y) \in \bigcup\limits_{m_1,m_2 \in
\mathbb{Z}}\, D_k+m_1+\imath m_2.
\end{equation}

The boundary conditions (\ref{bound_1})-(\ref{bound_2}) take the
form
\begin{equation}\label{bound_1a}
u_y|_{\partial Q^{(top)}_{(m_1,m_2)}}=-A
\sin \theta+q^+_1(x+m_1,m_1,m_2),
\end{equation}
\[
u_y|_{\partial Q^{(bottom)}_{(m_1,m_2)}}=-A
\sin \theta+q^-_1(x+m_1,m_1,m_2),
\]
and
\begin{equation} \label{bound_2a}
u_x|_{\partial Q^{(left)}_{(m_1,m_2)}}=-A
\cos \theta+q^-_2(y+m_2,m_1,m_2),
\end{equation}
\[
u_x|_{\partial Q^{(right)}_{(m_1,m_2)}}=-A
\cos \theta+q^+_2(y+m_2,m_1,m_2).
\]

The transmission conditions (\ref{bound_3}) and (\ref{bound_4}) along the inclusion interfaces can be
rewritten as follows:
\begin{equation}\label{bound_3a}
u=F_k(u_k),
\end{equation}
\begin{equation}\label{bound_4a}
\frac{\partial u}{\partial
n}=\frac{\partial u_k}{\partial n},
 \quad (x,y) \in \bigcup\limits_{m_1,m_2 \in \mathbb{Z}}\,(\partial D_k +m_1+\imath
 m_2),
\end{equation}
where the functions
\begin{equation}
\label{new_F}
F_k(\xi):=f(f_k^{-1}(\xi))
\end{equation}
are defined for all $\xi\in \mathbb{R}$.

Note that generally speaking the newly introduced function $u$ and $u_k$ are not continuous across the interface.

Zero mean value condition for the flux (\ref{bound_4aa}) can be
rewritten in the form:
\begin{equation}\label{bound_4aaa}
\int\limits_{\partial\, Q_{(m_1,m_2)}}\hspace{-5mm}\frac{\partial
u}{\partial n}ds=0.
\end{equation}
Moreover, since there is no source (sink)
inside the composite (neither in the matrix nor in any inclusion), the total heat flux through any closed
simply connected curve is equal to zero. Thus,
\begin{equation}\label{bound_4aak}
\int\limits_{\partial D_k +m_1+\imath
 m_2}\hspace{-5mm}\frac{\partial u_k}{\partial n}ds=0.
\end{equation}

As it follows from (\ref{bound_4a}), the same condition is valid for the function $u$ along
the boundary of each inclusion
\begin{equation}\label{bound_4aa0}
\int\limits_{\partial D_k +m_1+\imath
 m_2}\hspace{-5mm}\frac{\partial u}{\partial n}ds=0.
\end{equation}

Note that the function $F$ is also monotonic as it follows from the aforementioned arguments. Its derivative can be computed as follows
\begin{equation}
F'_k(\xi)=\frac{f'(f^{-1}_k(\xi))}{f'_k(f^{-1}_k(\xi))}=\frac{\lambda(T_k)}{\lambda_k(T_k)},
\end{equation}
where $\xi=f_k(T_k)$, and such representation is unique as the function $f_k$ is monotonic too.

Now we use the basic assumption of the paper on the nonlinear conduction coefficients
\begin{equation}\label{main_cond}
\lambda(T)= C_k\lambda_k(T).
\end{equation}
This property satisfies for any $T\in \mathbb{R}$ with some positive real constants $C_k$.
Then, one can immediately conclude that all functions $F_k$ are linear:
\begin{equation}\label{lineariz}
F_k(\xi)=D_k+C_k\xi.
\end{equation}
Note that from (\ref{Baiocchi}) we have $f(0)=0$ and $f_k(0)=0$, and, therefore, $D_k=0$.

Let us introduce inside the inclusions new harmonic functions:
\begin{equation}\label{u_tilda}
\tilde u_k(x,y)=C_ku_k(x,y).
\end{equation}
Then the transmission conditions (\ref{bound_3a}) and (\ref{bound_4a}) become
\begin{equation}\label{bound_3b}
u=\tilde u_k,
\end{equation}
\begin{equation}\label{bound_4b}
\frac{\partial u}{\partial
n}=\frac{1}{C_k}\frac{\partial \tilde u_k}{\partial n},
 \quad (x,y) \in \bigcup\limits_{m_1,m_2 \in \mathbb{Z}}\,(\partial D_k +m_1+\imath
 m_2),
\end{equation}

Thus, it is shown that the following statement is true.
\begin{theorem}
Let the assumption (\ref{main_cond}) be satisfied, then the nonlinear boundary value problem (\ref{1a})-(\ref{bound_2}), (\ref{bound_3}), (\ref{bound_4}) and
the linear boundary value problem (\ref{Lapl_a})-(\ref{bound_2a}), (\ref{bound_3b}), (\ref{bound_4b}) are equivalent.
\end{theorem}

In \cite{KPM} it has been shown that the problem (\ref{Lapl_a})-(\ref{bound_2a}), (\ref{bound_3b}), (\ref{bound_4b}) under the natural assumption
\begin{equation}
\label{super_condition}
q^\pm_j(t)=q^\pm_j(t+m_j,m_1,m_2),\quad
q^+_j(t)=q^-_j(t),\quad -1/2<t<1/2,
\end{equation}
is well-posed, and its solution $u$ possess the property that $\nabla u$ is a doubly periodic. This fact implies that the flux of the nonlinear boundary value problem (\ref{1a})-(\ref{bound_2}), (\ref{bound_3}), (\ref{bound_4}) is also doubly periodic provided the condition \eqref{super_condition} is fulfilled.

A new improved algorithm for solution of linear boundary value problem (\ref{Lapl_a})-(\ref{bound_2a}), (\ref{bound_3b})-(\ref{super_condition}) is developed and described in details in \cite{KPM}. We will use this approach in our computations.

\section{Average properties}

This section is devoted to evaluation of the average properties of
the nonlinear composite. We assume that the effective conductivity tensor $\Lambda_e$  depends on the average temperature
$\langle T \rangle$ and  is defined in the following way:
\begin{equation}
\label{average}
\langle\lambda(T)\nabla T\rangle= \Lambda_e(\langle T \rangle)\langle\nabla
T\rangle,\quad \mbox{or}\quad
R_e(\langle T \rangle)\langle\lambda(T)\nabla T\rangle=\langle\nabla T\rangle,
\end{equation}
where $\Lambda_{e}$ is the tensor of the composite conductivity while $R_e=\Lambda_{e}^{-1}$ is the effective resistance tensor (for similar definition see \cite{SZ}). Here, the operator operator $\langle \cdot \rangle$ is defined as the integral over the cell volume in the standard form:
\[
\langle f \rangle=\iint\limits_{Q_{(m_1,m_2)}}\,f(x,y)\,dxdy.
\]
Note that such definition needs further justification
as question arises whether the approach is invariant with respect to the averaging cell. We will discuss this issue later during the computations.

Let
us compute all the terms involving the constructed solution. The total flux in $x$-direction can be transformed to
\[
 \iint\limits_{Q_{(m_1,m_2)}}\,\lambda(T)\frac{\partial
T}{\partial x}\,dxdy =\iint\limits_{D_0+m_1+\imath m_2}\,\lambda(T)\frac{\partial
T}{\partial x}\,dxdy +\sum\limits_{k=1}^N\,
\iint\limits_{D_k+m_1+\imath m_2}\lambda_k(T_k) \frac{\partial T_k}{\partial x}\,
dxdy
\]
\[
=\iint\limits_{Q_{(m_1,m_2)}}\,(f(T))_x\,dxdy+ \sum\limits_{k=1}^N\,
\iint\limits_{D_k+m_1+\imath m_2} \, (f_k(T_k))_x\,
dxdy
\]
\begin{equation}\label{lam_x1}
=\iint\limits_{D_0+m_1+\imath m_2}\,\frac{\partial u}{\partial x}\,dxdy+
\sum\limits_{k=1}^N\, \iint\limits_{D_k+m_1+\imath m_2} \, \frac{\partial
u_k}{\partial x}\, dxdy.
\end{equation}
Using the first Green's formula and formulas (\ref{Lapl_a}), (\ref{Lapl_b}) and (\ref{bound_4a}), we obtain
$$
\iint\limits_{Q_{(m_1,m_2)}}\,\lambda(T)\frac{\partial
T}{\partial x}\,dxdy =-A\cos\theta,
$$
and
$$
\iint\limits_{Q_{(m_1,m_2)}}\,\lambda(T)\frac{\partial
T}{\partial y}\,dxdy=-A\sin\theta.
$$
Note that the last two identities are trivial consequences of the assumption on the heat flux in the composite in the absence of the sources and/or sinks.

Thus,
\begin{equation}
\label{ave_flux} \langle\lambda(T)\,\nabla
T\rangle=-A[\cos\theta,\sin\theta]^\top.
\end{equation}

Taking into account Gauss-Ostrogradsky formula and (\ref{bound_3}), the components of the term $\langle\nabla T\rangle$ in (\ref{average}) are defined as
$$
\iint\limits_{Q_{(m_1,m_2)}}\,\frac{\partial T}{\partial
x}\,dxdy
=\iint\limits_{D_0+m_1+\imath m_2}\,\frac{\partial T}{\partial
x}\,dxdy+\sum\limits_{k=1}^N\, \iint\limits_{D_k+m_1+\imath m_2}\,\frac{\partial T_k}{\partial
x}\,dxdy
$$
$$
=\oint\limits_{\partial D_0+m_1+\imath m_2} T(s) \cos(n_s,e_i) \,ds+\sum\limits_{k=1}^N\, \oint\limits_{\partial D_k+m_1+\imath m_2} [T_k(s)-T(s)] \cos(n_s^k,e_i) \,ds
$$
$$
=\oint\limits_{\partial D_0+m_1+\imath m_2} T(s) \cos(n_s,e_i) \,ds=\oint\limits_{\partial D_0+m_1+\imath m_2} \,f^{-1}(u(x,y)) \cos(n_s,e_i) \,ds,
$$
where $n_s$ and $n_s^k$ are the outward unit normal vectors to $\partial D_0+m_1+\imath m_2$ and $\partial D_k+m_1+\imath m_2$, respectively, and $e_i$ is the basis vector.
Analogously,
$$
\iint\limits_{Q_{(m_1,m_2)}}\,\frac{\partial T}{\partial
y}\,dxdy
=\iint\limits_{D_0+m_1+\imath m_2}\,\frac{\partial T}{\partial
y}\,dxdy+\sum\limits_{k=1}^N\, \iint\limits_{D_k+m_1+\imath m_2}\,\frac{\partial T_k}{\partial
y}\,dxdy
$$
$$
=\oint\limits_{\partial D_0+m_1+\imath m_2} T(s) \cos(n_s,e_j) \,ds=\oint\limits_{\partial D_0+m_1+\imath m_2} \,f^{-1}(u(x,y)) \cos(n_s,e_j) \,ds.
$$

Finally, the average temperature is
\begin{multline}\label{av_T}
\langle T \rangle=\iint\limits_{Q_{(m_1,m_2)}}\,T(x,y)\,dxdy \\
=\iint\limits_{D_0+m_1+\imath m_2}\,f^{-1}(u(x,y))\,dxdy+
\sum_{k=1}^N\,\iint\limits_{D_k+m_1+\imath m_2}\,f^{-1}_k(u_k(x,y))\,dxdy.
\end{multline}

Bearing in mind these computations, it is more convenient to compute  the components of the flux resistance tensor $R_e$ from (\ref{average}) first and than
to determine the effective conductivity tensor $\Lambda_{e}=R_e^{-1}$.


\section{Numerical example}

\subsection{Description of the periodic composite}

We consider a composite where four inclusions are situated inside the cell $Q_{(0,0)}$ with the centers: $a_1=-0.18+0.2 \imath$, $a_2=0.33-0.34 \imath$, $a_3=0.33+0.35 \imath$, $a_4=-0.18-0.2\imath$. For calculation  we assume that radii are the same $r_k=R=0.145$ (inclusions are very close to each other in neighboring cells). Note that the volume fraction of the inclusions for such composite is $0.2642$.
Further, we take the following characteristics of the flux: $\theta=0, A=-1$, thus the heat flows in $x$-direction.

\begin{figure}[h!]
\begin{center}
\resizebox*{7cm}{!}{\includegraphics{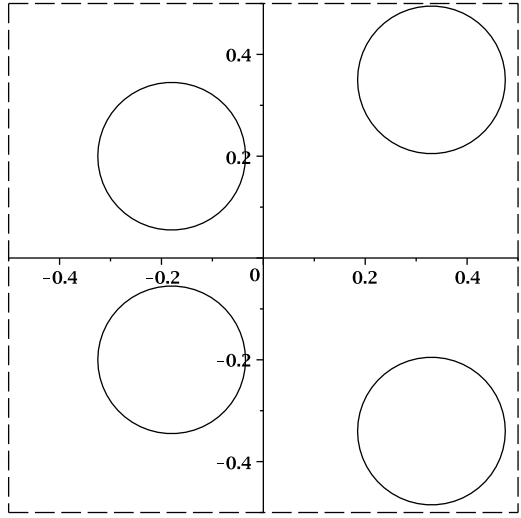}}
\caption{Configuration of the unit cell with four inclusions considered in computation}%
\label{fig:4_inclusions}
\end{center}
\end{figure}

Computation will be done with the use of the algorithm described in \cite{KPM}, where we
sought a solution in the form of Taylor series. For the chosen configuration, we take $M=6$ first items in the Taylor series which guarantees
the accuracy for the temperature with the error less than $10^{-6}$.

We choose the conductivities $\lambda(T)$ and $\lambda_k(T)$ in such a way that the condition (\ref{main_cond}) is satisfied.
Namely,  we consider the conductivities given in the following form
\begin{equation}\label{lam_ex}
\lambda(T)=
\begin{cases}
y_1,\quad T<x_1, \cr
y_2+\frac{y_1-y_2}{x_1}\, T,\quad x_1\leq T \leq 0, \cr
y_2+\frac{y_1-y_2}{x_2}\, T,\quad 0\leq T \leq x_2, \cr
y_1,\quad T> x_2,
\end{cases}
\end{equation}
\begin{equation}\label{lam_k_ex}
\lambda_k(T)=
\begin{cases}
y_3,\quad T<x_1, \cr
y_4+\frac{y_3-y_4}{x_1}\, T,\quad x_1\leq T \leq 0, \cr
y_4+\frac{y_3-y_4}{x_2}\, T,\quad 0\leq T \leq x_2, \cr
y_3,\quad T> x_2,
\end{cases}
\end{equation}
where let $y_1,y_2,y_3,y_4$ be positive constants and $x_1<x_2$.

Then,
\begin{equation} \label{f_ex}
f(T)=
\begin{cases}
y_1 T+x_1 \frac{y_2-y_1}{2},\quad T<x_1, \cr
\frac{y_1-y_2}{2 x_1}T^2+y_2 T, \quad x_1\le T\le 0, \cr
\frac{y_1-y_2}{2 x_2}T^2+y_2 T,\quad 0\le T\le x_2, \cr
y_1 T+x_2 \frac{y_2-y_1}{2}, \quad T>x_2.
\end{cases}
\end{equation}
The function $f_k$ has the same form. For $f_k$, we take $y_3$ instead of $y_1$ and $y_4$ instead of $y_2$.
Calculating $f_k^{-1}$, we obtain
\begin{equation} \label{invf_k_ex}
f_k^{-1}(\xi)=
\begin{cases}
\frac{\xi-\frac{x_1(y_4-y_3)}{2}}{y_3},\quad \xi<\frac{x_1(y_3+y_4)}{2}, \cr
\frac{-x_1 y_4-\sqrt{(x_1 y_4)^2+2(y_3-y_4)x_1 \xi}}{y_3-y_4},\quad \frac{x_1(y_3+y_4)}{2}\le \xi\le 0, \\[1mm]
\frac{-x_2 y_4+\sqrt{(x_2 y_4)^2+2(y_3-y_4)x_2 \xi}}{y_3-y_4},\quad 0\le \xi \le \frac{x_2(y_3+y_4)}{2}, \cr
\frac{\xi-\frac{x_2(y_4-y_3)}{2}}{y_3},\quad \xi >\frac{x_2(y_3+y_4)}{2}.
\end{cases}
\end{equation}
Note that all required properties hold true: $\lambda(T)$ is continuous, $f'(T)=\lambda(T)$ and $f^{-1}(f(x))=x$.

\begin{figure}[h!]
\begin{center}
\resizebox*{8cm}{!}{\includegraphics{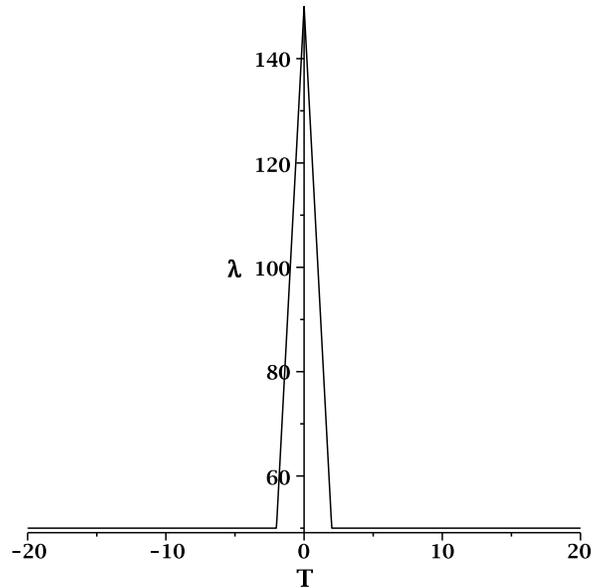}}
\caption{The function $\lambda_k$.}%
\label{fig:cond}
\end{center}
\end{figure}

If $y_2=s y_1$ and $y_4=s y_3$ with an arbitrary $s \in \mathbb{R}$ then the function $F$ defined as $F(\xi_k):=f(f_k^{-1}(\xi_k))$ has the following form
\begin{equation}\label{F_k_lin}
F(\xi_k)=\frac{y_1}{y_3}\, \xi_k, \qquad  -\infty<x_1<x_2<+\infty
\end{equation}
with $C_k=y_1/y_3$. We take for the calculations $x_1=-2,\, x_2=2$, and assume that
 $y_1=4.5,\, y_2=13.5, \, y_3=50, \, y_4=150$, i.e., $s=3$ and $C_k=0.09$. On Fig. \ref{fig:cond}, we represent the conductivity function $\lambda_k$ of the inclusions.
The function $\lambda$ has the identical shape with the pike taking value $\lambda(0)=0.09 \cdot \lambda_k(0)=13.5$.

\subsection{Evaluation of the effective conductivity tensor}

Note that in the linear case the temperature is defined up to an arbitrary additive constant (see \cite{KPM}), and
this constant is not involved in the determination of the effective effective conductivity of the composite material.
In the nonlinear case this is not generally speaking the case, and one needs to clarify how the additive constant  appears on the stage of solving
 the auxiliary linear problem (\ref{Lapl_a})-(\ref{bound_2a}), (\ref{bound_3b})-(\ref{super_condition}) with respect to the functions $u$ and $u_k$ influencing (or not) on the computations of effective conductivity tensor of the equivalent nonlinear composite.

Two procedures can be suggested to evaluate the effective conductivity.

\begin{itemize}
\item
First, one can solve the auxiliary linear boundary value problem in the doubly periodic domain preserving its uniqueness
by any appropriately chosen condition (for example, here we assume that the function $u=u_*$ satisfies the condition $u_*(0)=0$). Then, to
evaluate the properties of the composite material, one can compute the average temperature and the average resistivity for each particular unit cell
presenting the data as the functional relationship $R_e=R_e(\langle T \rangle)$.

It is clear from the character of the chosen
conductivities that the domain where the nonlinear behavior manifests itself lies only inside an infinite strip of unknown \emph{finite thickness} which depends on the flux intensity $A$ and its direction.
Thus, it is not a surprise that the effective conductivity tensor demonstrates nonlinear behavior within a finite interval of temperatures, and thus,
one does not need to trace all the cells.
On the other hand, there is still infinite number of the cells belonging to the strip, and therefore one can expect that the result of such procedure is representative enough.

\item
Another method for evaluation of the average properties consists in the following. One can consider an arbitrary cell in the original domain and build a set of solutions of the auxiliary problem in the form $u=u_*+C$, where $C$ is an arbitrary constant. Then, for every constant $C$, the components of the resistance tensor and the average of the temperature computed in Section 4 are functions of the parameter $C$. Changing them continuously from $-\infty$ to $\infty$, one receives the sought effective conductivity tensor of the composite as a continuous function of the effective temperature. Moreover, for the conductivities of the composite components analyzed in this example, the nonlinear character of the relationship will be observed only within the finite interval of the parameter $C$. One can realized that this procedure does not depend on the chosen cell.
\end{itemize}

Note that the both methods allow one to determine two components of the resistance tensor $R_e$ for each perpendicular flux direction.
Thus considering $\theta=0$, we define $R_e[1,j]=R_e[1,j](\langle T\rangle)$ ($j=1,2$), while choosing $\theta=\pi/2$ we find $R_e[2,j]=R_e[2,j](\langle T\rangle)$. As a result, the entire tensor $R_e(\langle T\rangle)$ is defined.

\begin{figure}[h!]
\begin{center}
\resizebox*{13cm}{!}{{\includegraphics{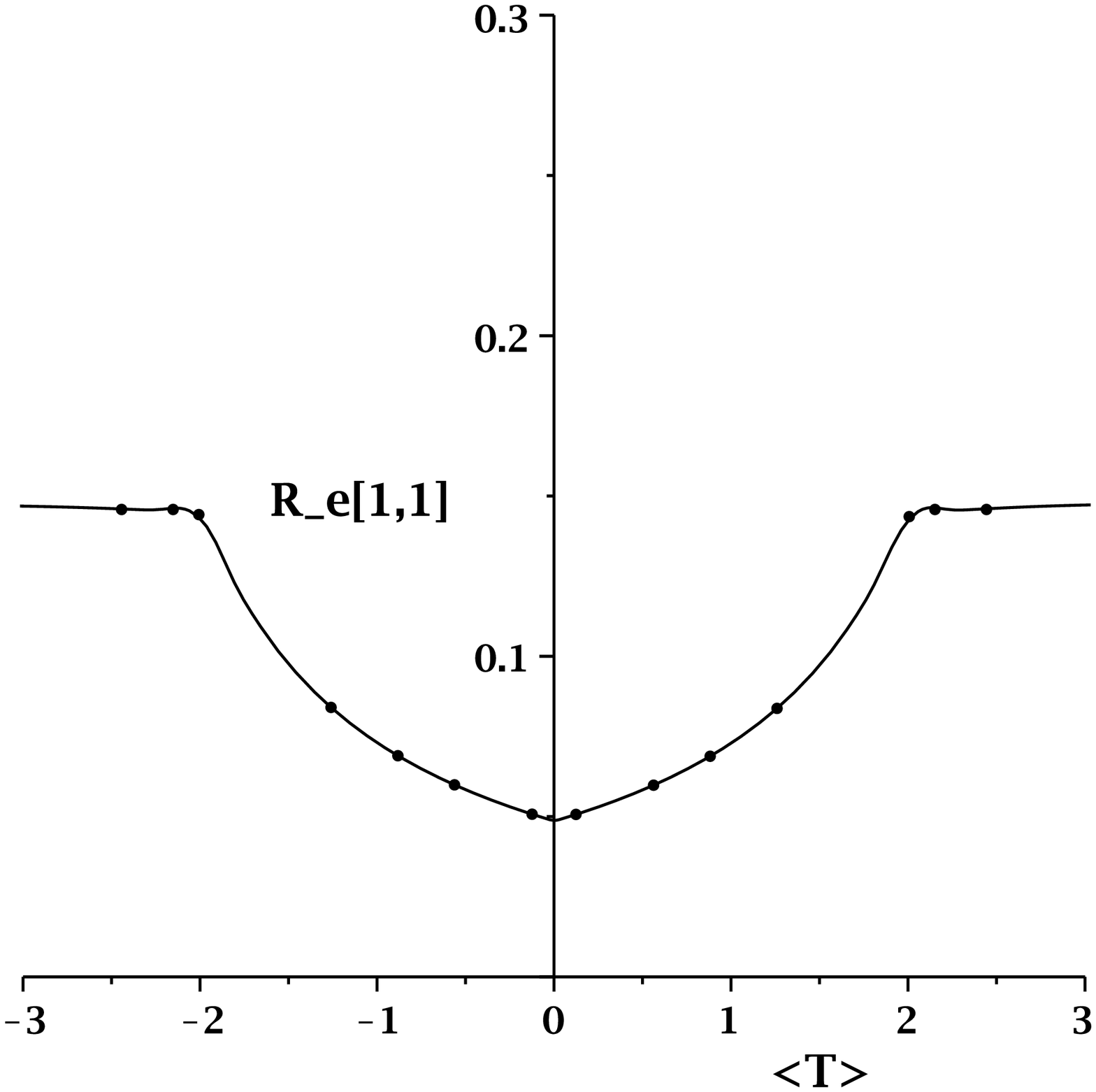},\includegraphics{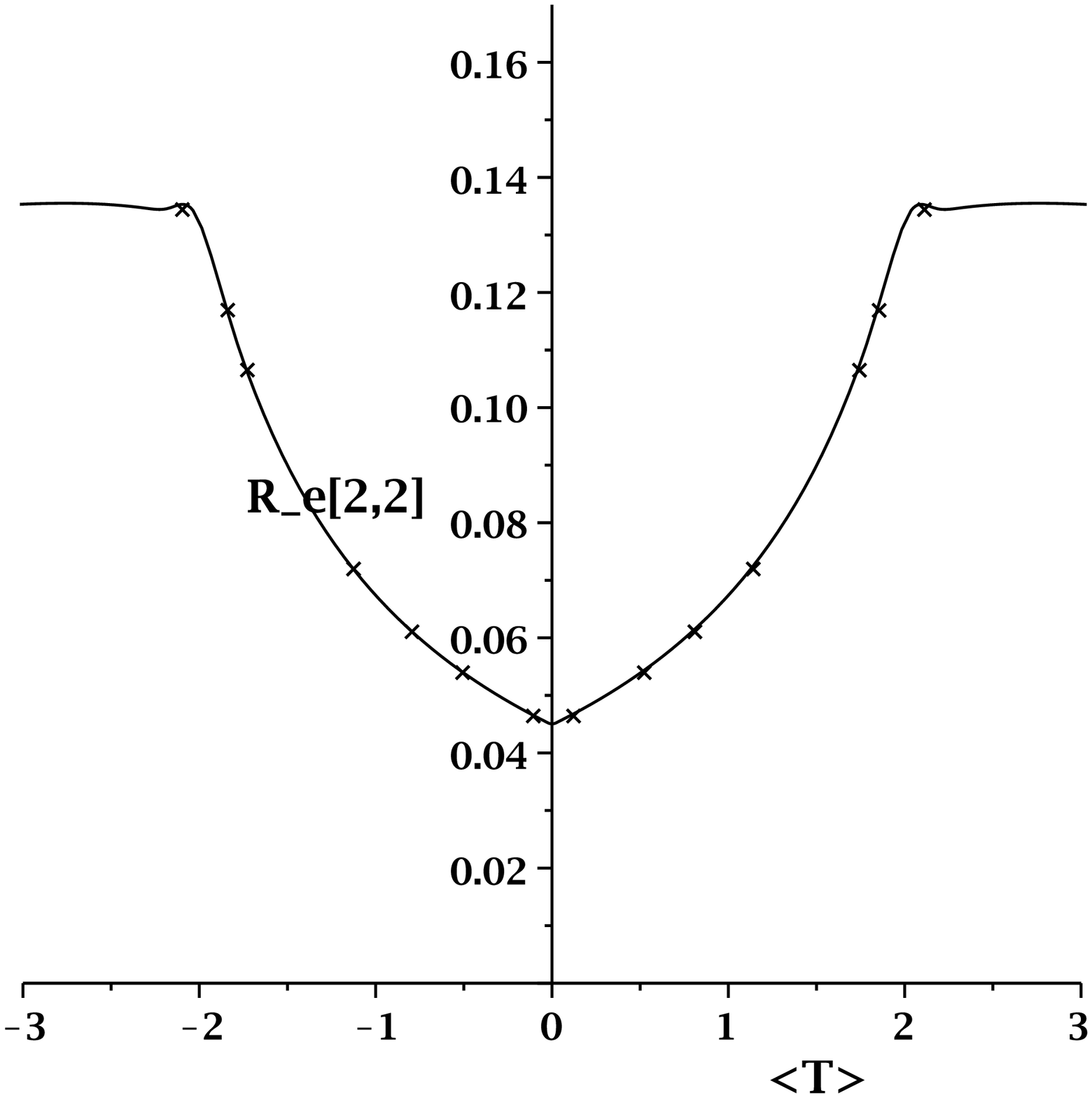}}}
\caption{Main diagonal elements of the average resistance tensor $R_e$ ($R_e[1,1]$ and $R_e[2,2]$) computed by each of the proposed methods. Dots corresponds to the
values computed in different cells.}%
\label{fig:R_1122}
\end{center}
\end{figure}

\begin{figure}[h!]
\begin{center}
\resizebox*{13cm}{!}{{\includegraphics{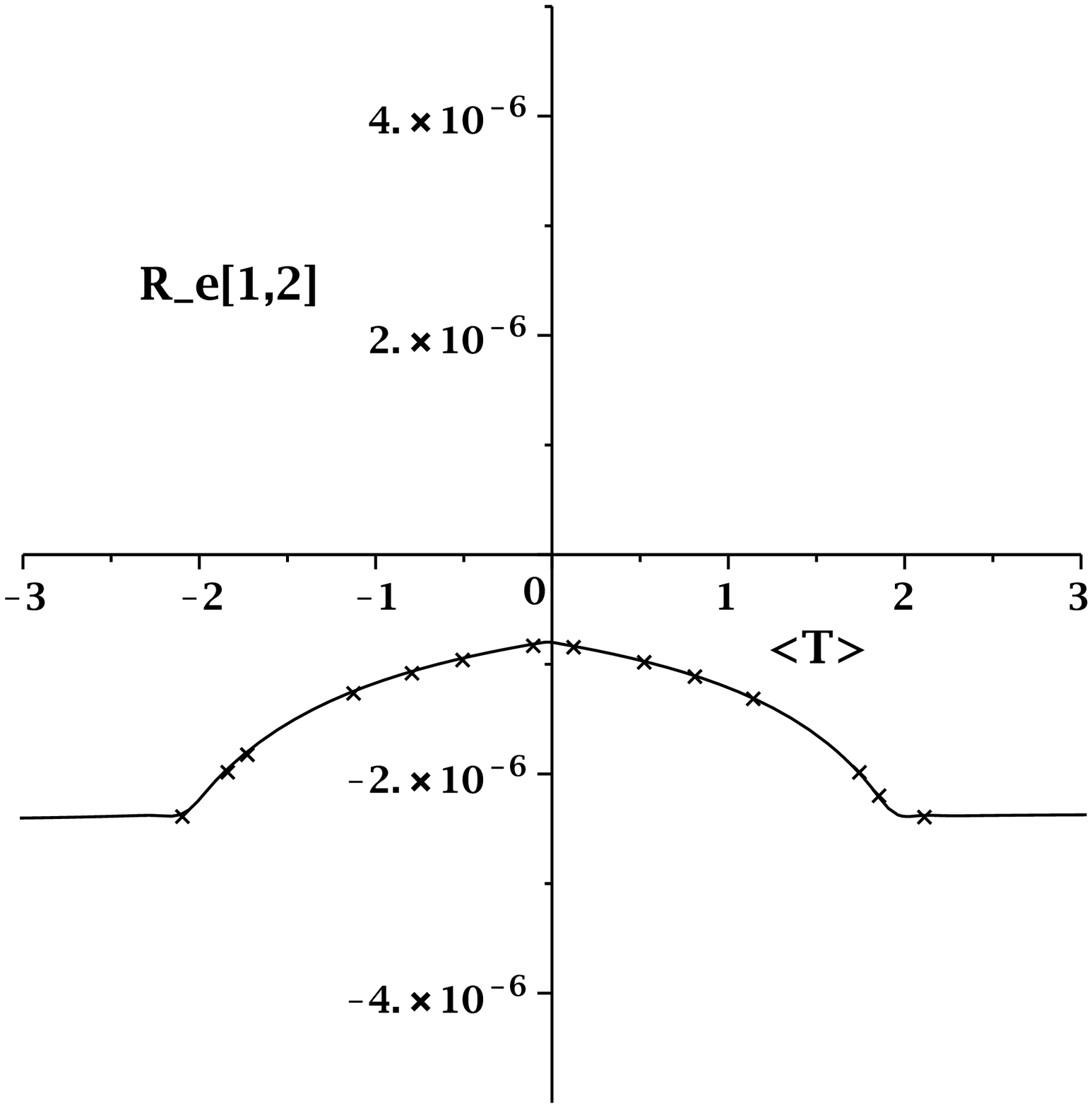},\includegraphics{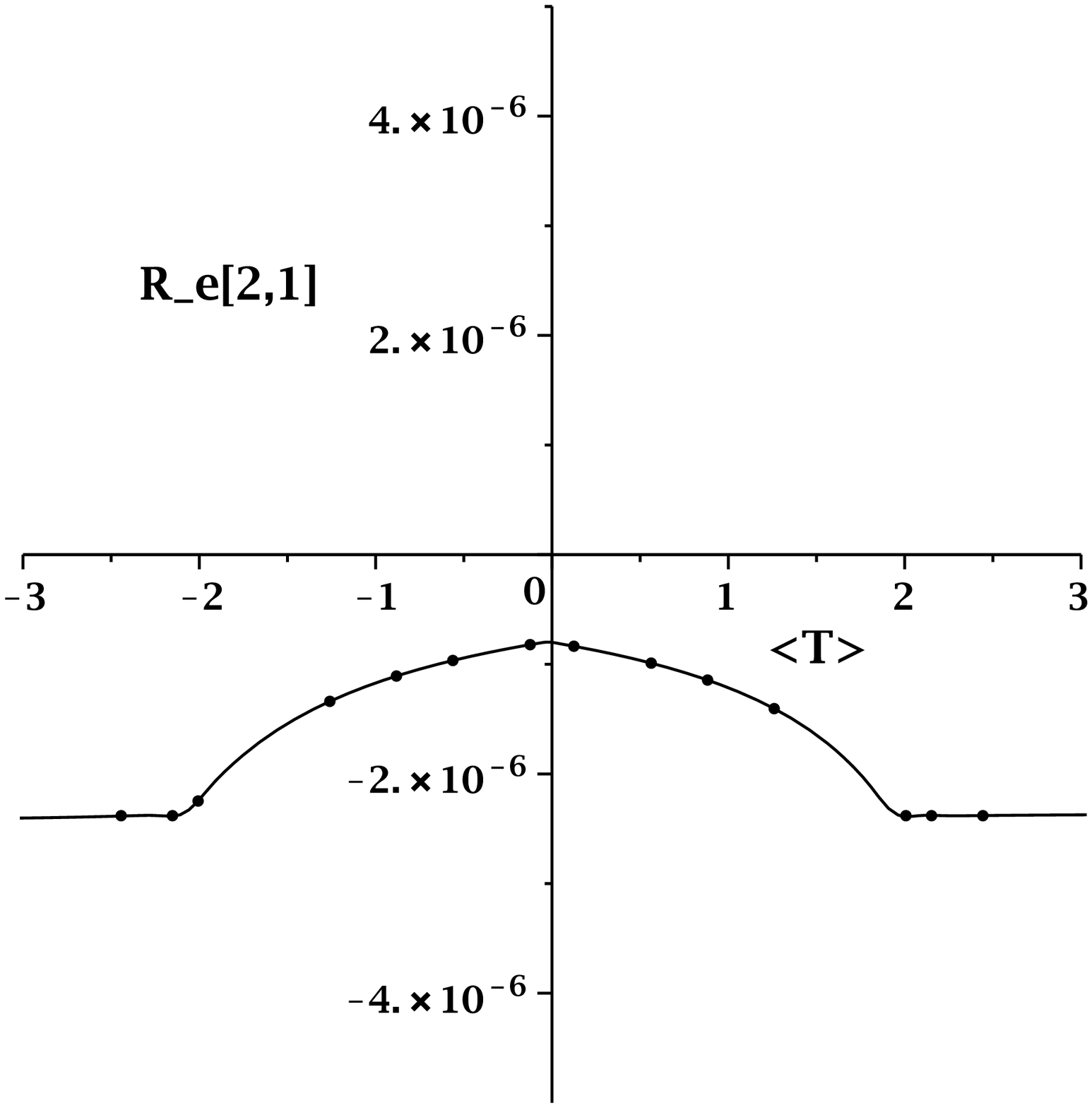}}}
\caption{Components $R_e[1,2]$ and $R_e[2,1]$ of the average resistance tensor $R_e$ computed by each of the proposed methods.}
\label{fig:R_1221}
\end{center}
\end{figure}

To demonstrate that these two aforementioned procedures are equivalent, we use both of them in our computations. The respective
components of the effective resistance tensor are represented on Figures \ref{fig:R_1122},\ref{fig:R_1221}.
Dots on the curves correspond to the second approach, while the continuous lines correspond to values computed for consecutive unit cells.
These continuous lines were obtained by spline interpolations. One can expect that, due to the chosen functions determining the conductivities of the components, the effective conductivity should be an even function of the average temperature. However, there is no reason to restrict ourselves to the computations for the cells where the average temperature is positive (negative) as this decrease the amount of information allowing to draw the curve more accurate.

Discrepancy between the methods was on the level $10^{-5}$ while the computational accuracy of the solution itself was $10^{-6}$. Taking into account that fact that one needs to integrate and interpolate the data to compute the average properties, this can be considered as a good evidence that the both methods provide the same results. However, an accurate mathematical proof is still to be delivered.

\newpage

\begin{figure}[h!]
\begin{center}
\resizebox*{13cm}{!}{{\includegraphics{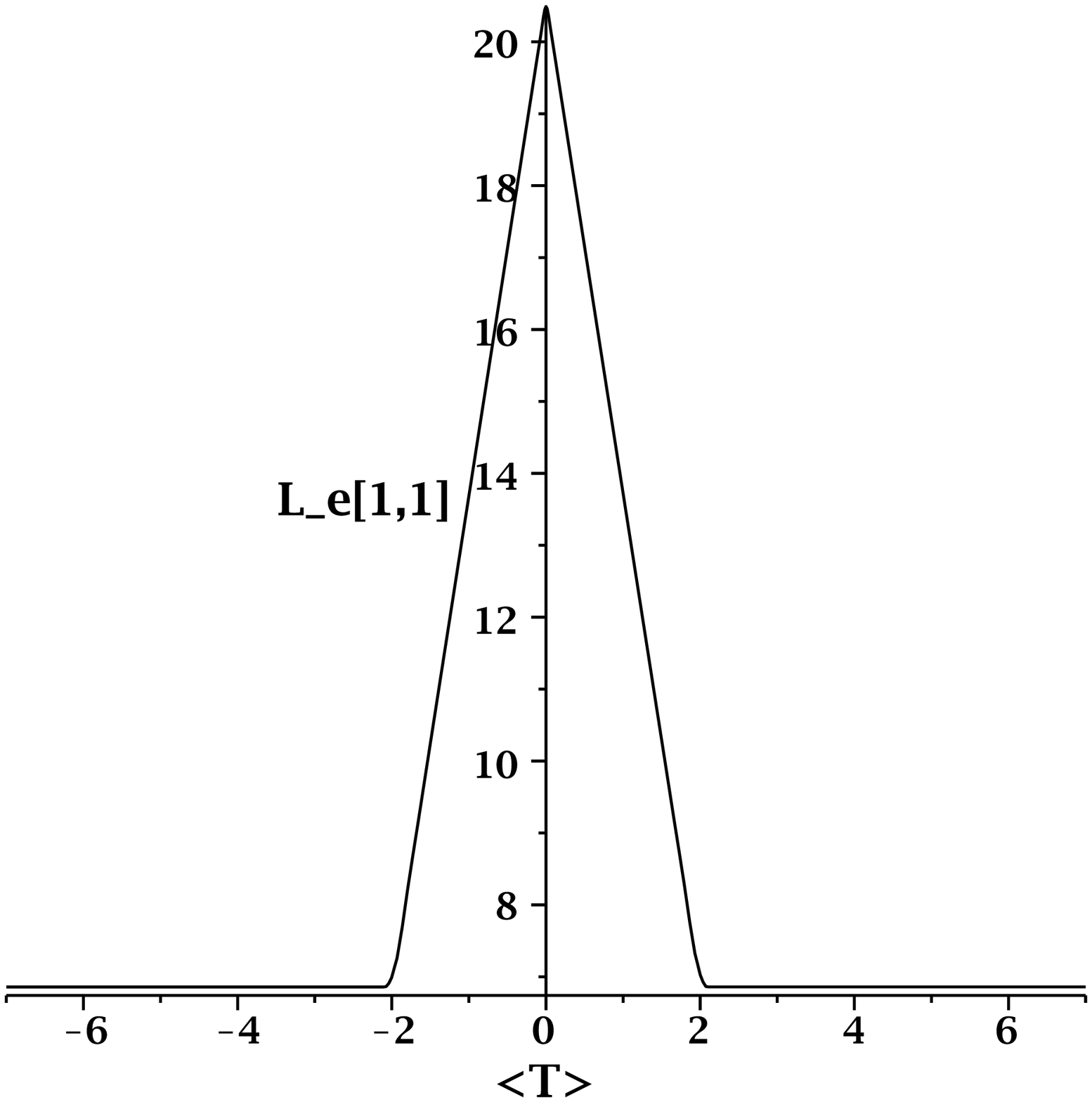},\includegraphics{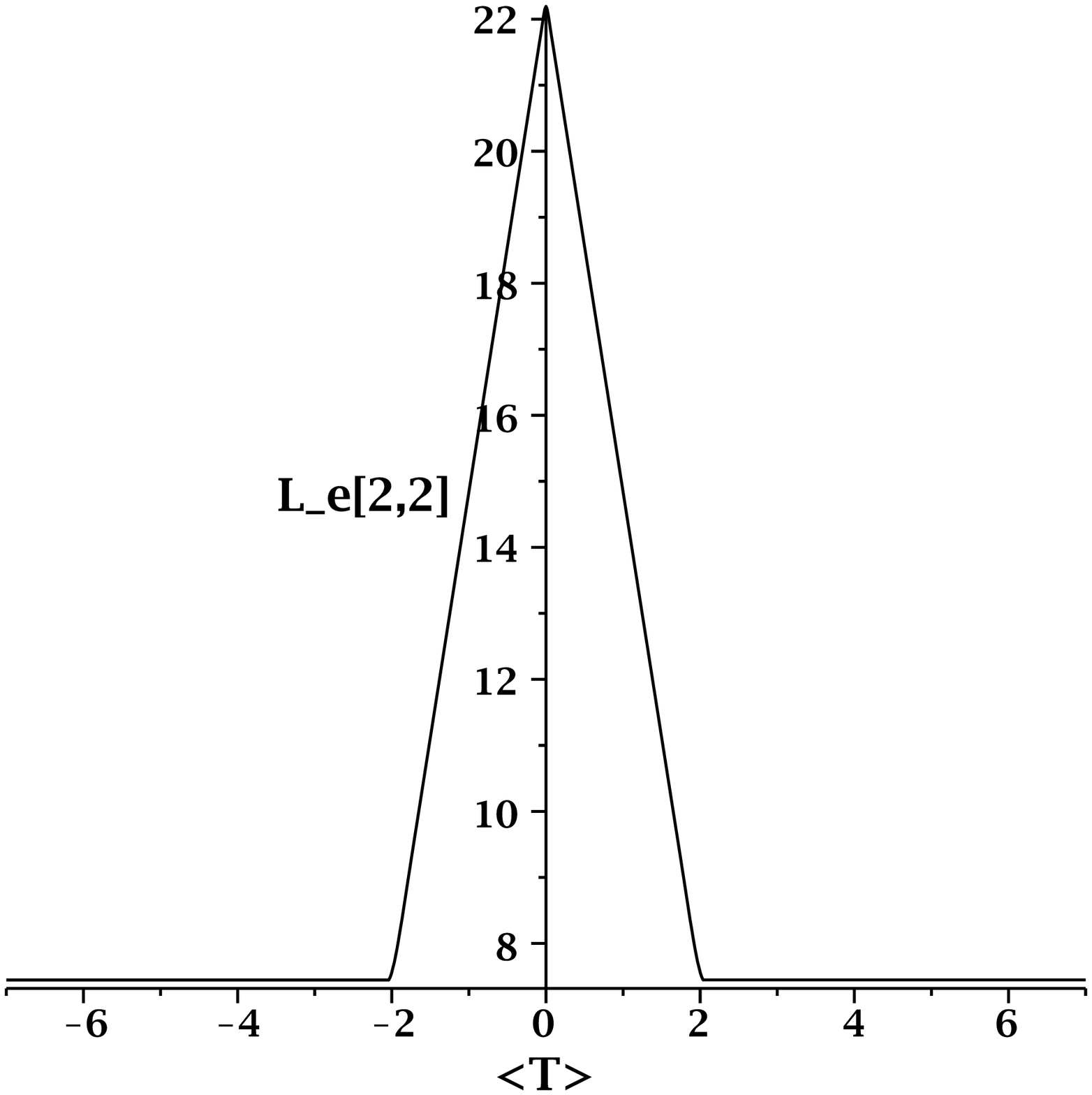}}}
\caption{Diagonal components $\Lambda_e[1,1]$, $\Lambda_e[2,2]$ of the effective conductivity tensor $\Lambda_e$ as function of the temperature average $\langle T\rangle$.}%
\label{fig:L_1122}
\end{center}
\end{figure}

\begin{figure}[h!]
\begin{center}
\resizebox*{7cm}{!}{\includegraphics{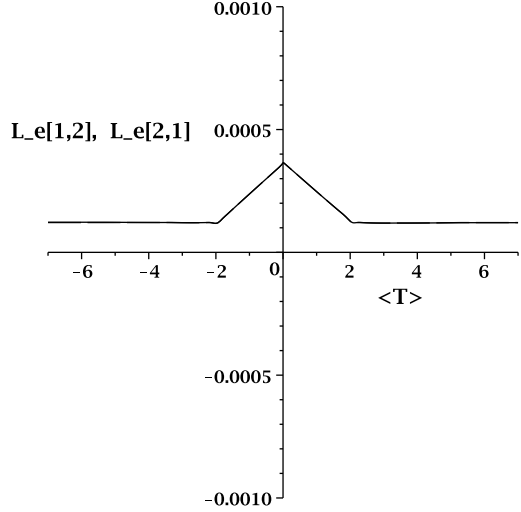}}
\caption{Components $\Lambda_e[1,2]$ and $\Lambda_e[2,1]$ of the effective conductivity tensor.}
\label{fig:L_1221}
\end{center}
\end{figure}

Finally, having the resistance tensor $R_e(\langle T\rangle)$, we calculate the effective conductivity tensor $\Lambda_e(\langle T\rangle)$ defined in (\ref{average})
as a matrix function of the average temperature $\langle T \rangle$. The respective results are represented on Figures \ref{fig:L_1122},\ref{fig:L_1221}, respectively. One can see that the shape of the functions are quite similar to that demonstrated by the function $\lambda(\langle T\rangle)$ and $\lambda_k(\langle T\rangle)$ for the composite components. Only some deviations
can be observed near the points when the functions are not smooth.


\newpage

\subsection{Hashin-Shtrikman bounds and other estimates}

In \cite{GalTelTok}, \cite{MiKo} the general estimates for the average composite properties have been evaluated in their more general nonlinear formulation.
In particular, the elementary bounds of the effective conductivity tensor in our notations has the form:
\begin{equation}\label{el_bounds}
\mu_1(T)I\leq \Lambda_e(T) \leq \mu_2(T)I,
\end{equation}
where $I$ is the unit tensor and
$$
\mu_1 (T)=\left(\frac{1-N \pi R^2}{\lambda(T)}+\frac{N \pi R^2}{\lambda_k(T)}\right)^{-1}= \left(\frac{0.7358}{\lambda(T)}+\frac{0.2642}{\lambda_k(T)}\right)^{-1},
$$
$$
\mu_2 (T)=\lambda(T)(1-N \pi R^2)+\lambda_k(T)N \pi R^2=0.7358\cdot \lambda(T)+0.2642\cdot\lambda_k(T).
$$
Inequalities (\ref{el_bounds}) are the so-called the Reuss-type and Voigt-type bounds on the
effective coefficients, cf. \cite{GalTelTok}, \cite{Kusch}.

Note that we have constructed an analytical solution to the nonlinear problem and directly determined the effective conductivity of the composite. However,
since the question on how the average properties of such composites should be determined and understood, we decided to verify our computations relate to the general
estimates following from the variational analysis. Although, to estimate (\ref{el_bounds}) is rather crude, we start from it.

Let us recall that if $A$ and $B$ are matrices,
then the notation $A\ge B$ means that inequality $(Ax,x)\ge(Bx,x)$ holds true for an arbitrary vector $x\in \mathbb{R}^n$ ($n=2$ in our case).
In other words one needs to show that the following inequalities are true:
$$
m_{11}(T)=\mu_1-\lambda^e_{11}\le0 , \quad m_{21}(T)=\mu_2-\lambda^e_{11}\ge0,
$$
$$
m_{12}(T)=4(\mu_1-\lambda^e_{11})(\mu_1-\lambda^{e}_{22})-(\lambda^{e}_{12}+\lambda^{e}_{21})^2\ge0,
$$
$$
m_{22}(T)=4(\mu_2-\lambda^e_{11})(\mu_2-\lambda^{e}_{22})-(\lambda^{e}_{12}+\lambda^{e}_{21})^2\ge0.
$$

\begin{figure}[h!]
\begin{center}
\resizebox*{14cm}{!}{{\includegraphics{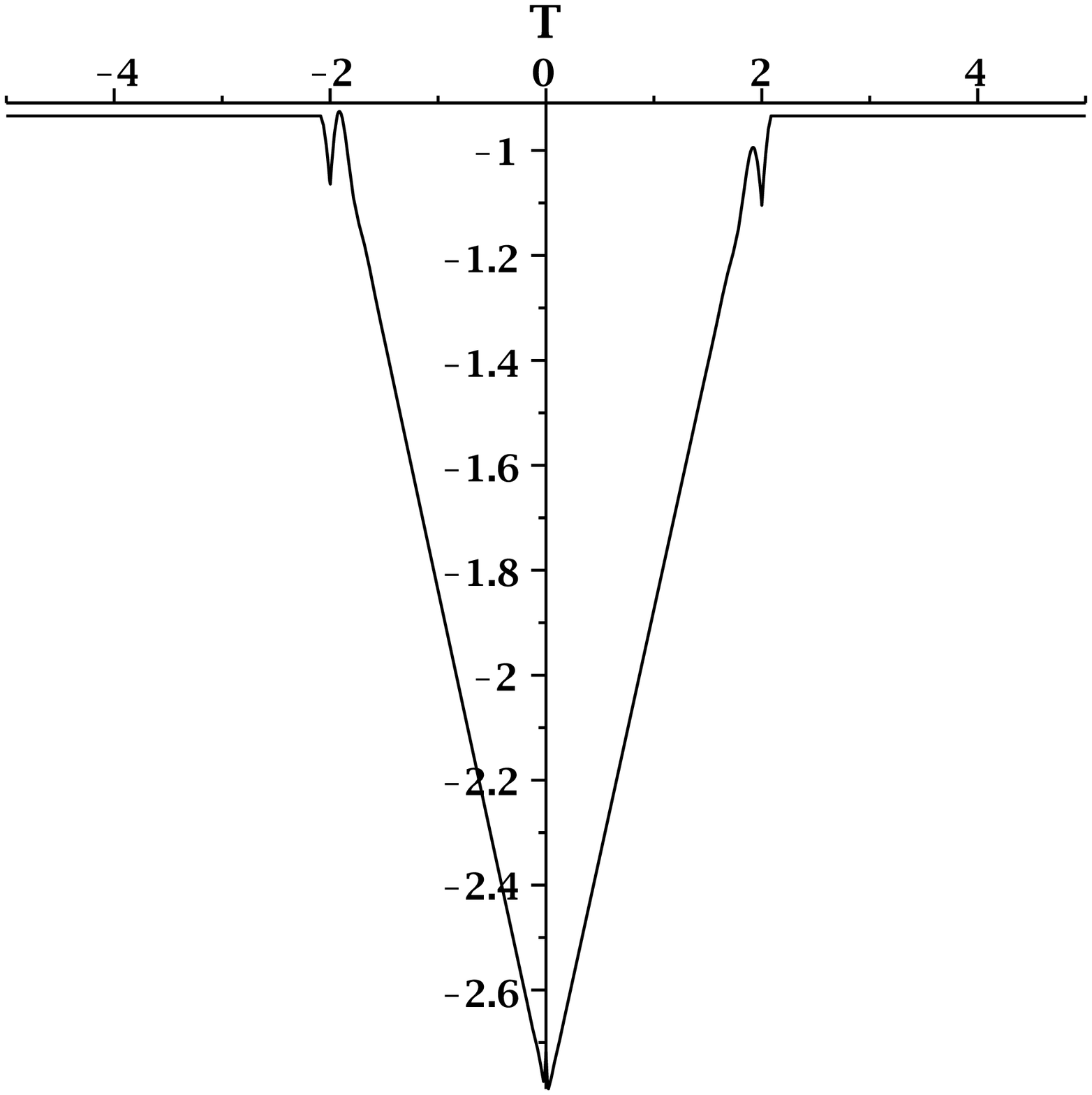},\includegraphics{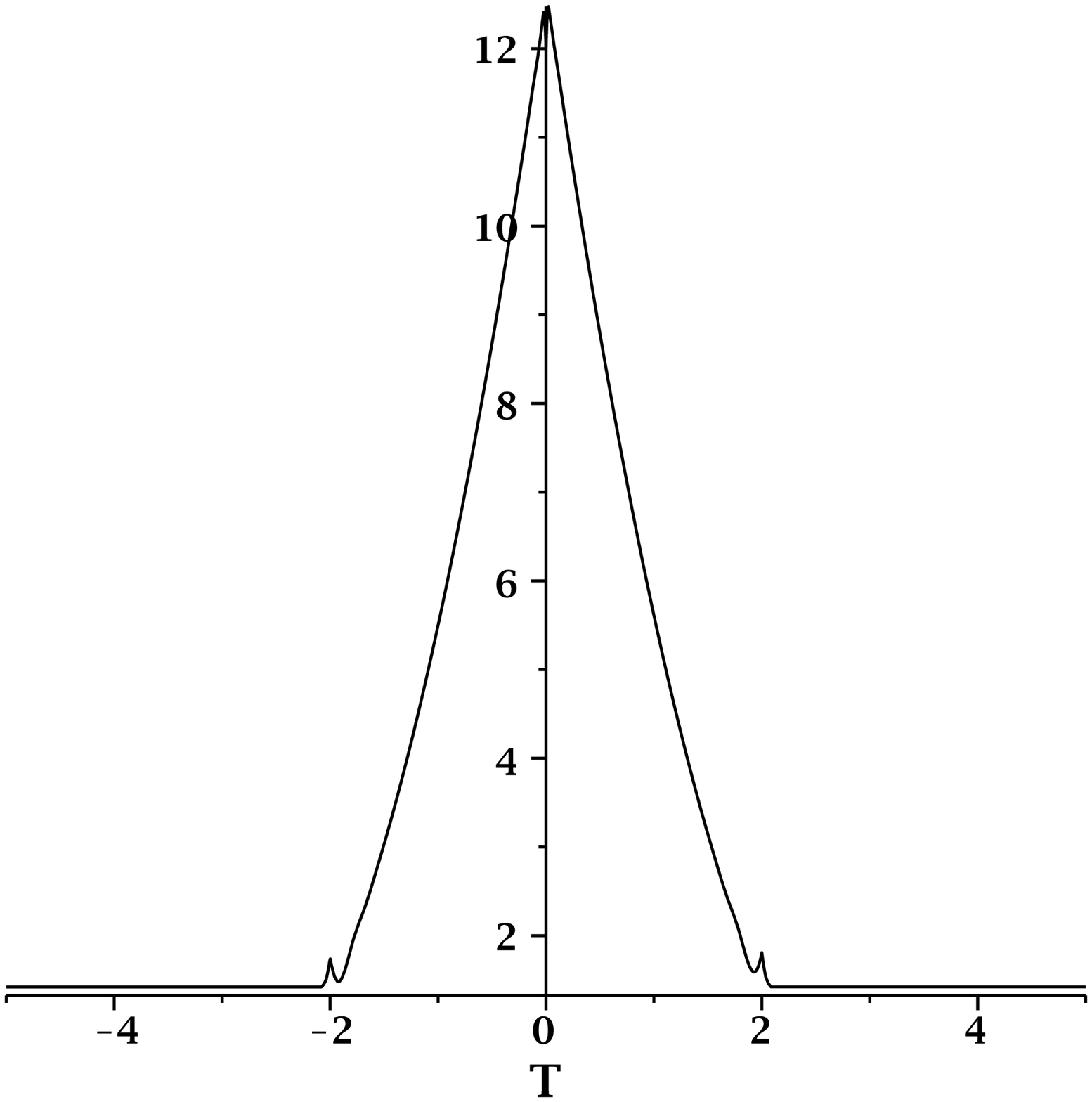}}}
\caption{Verification of the Reuss-type inequality for the effective properties of the nonlinear composite.}%
\label{minor1}
\end{center}
\end{figure}
\begin{figure}[h!]
\begin{center}
\resizebox*{14cm}{!}{{\includegraphics{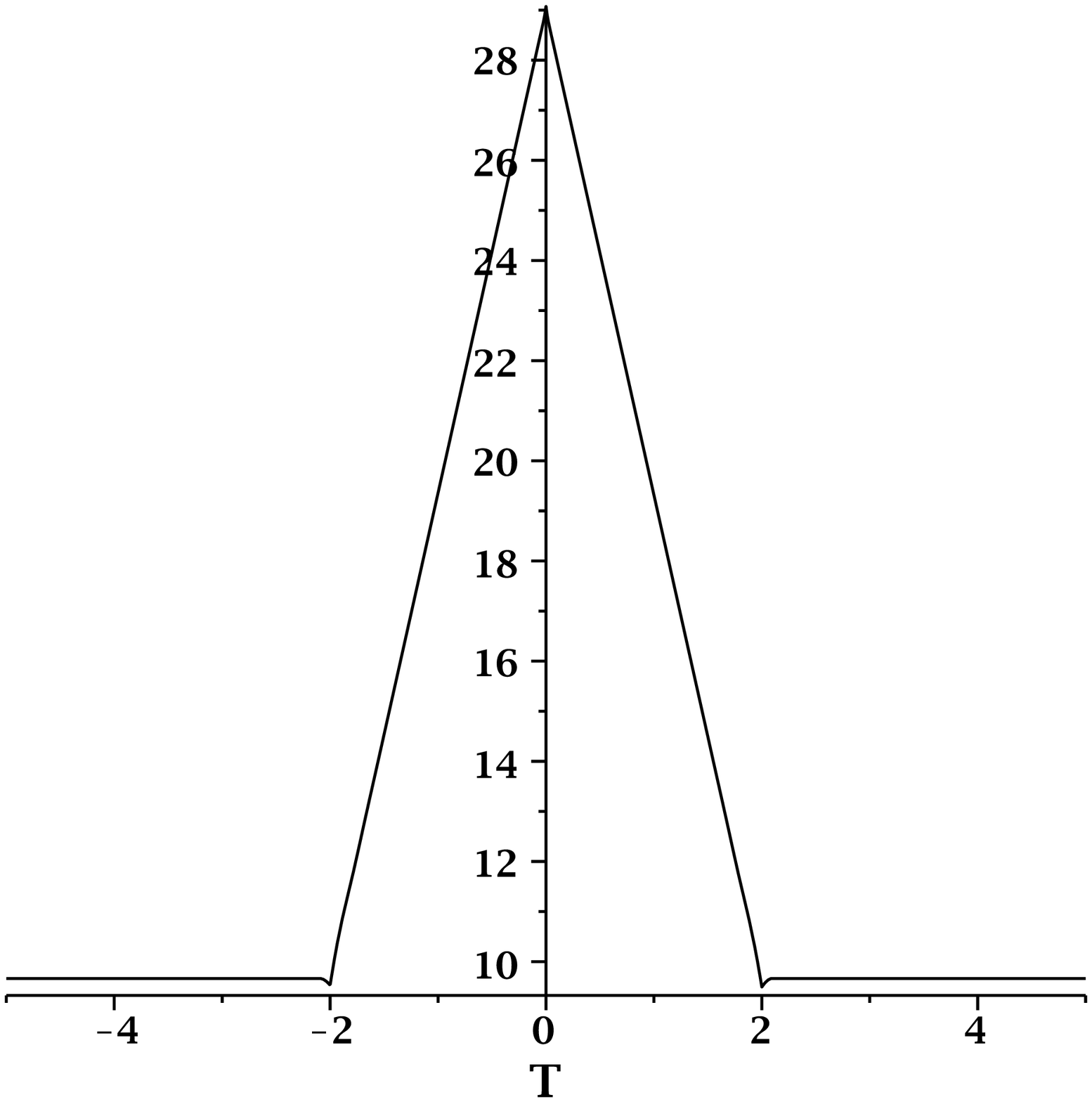},\includegraphics{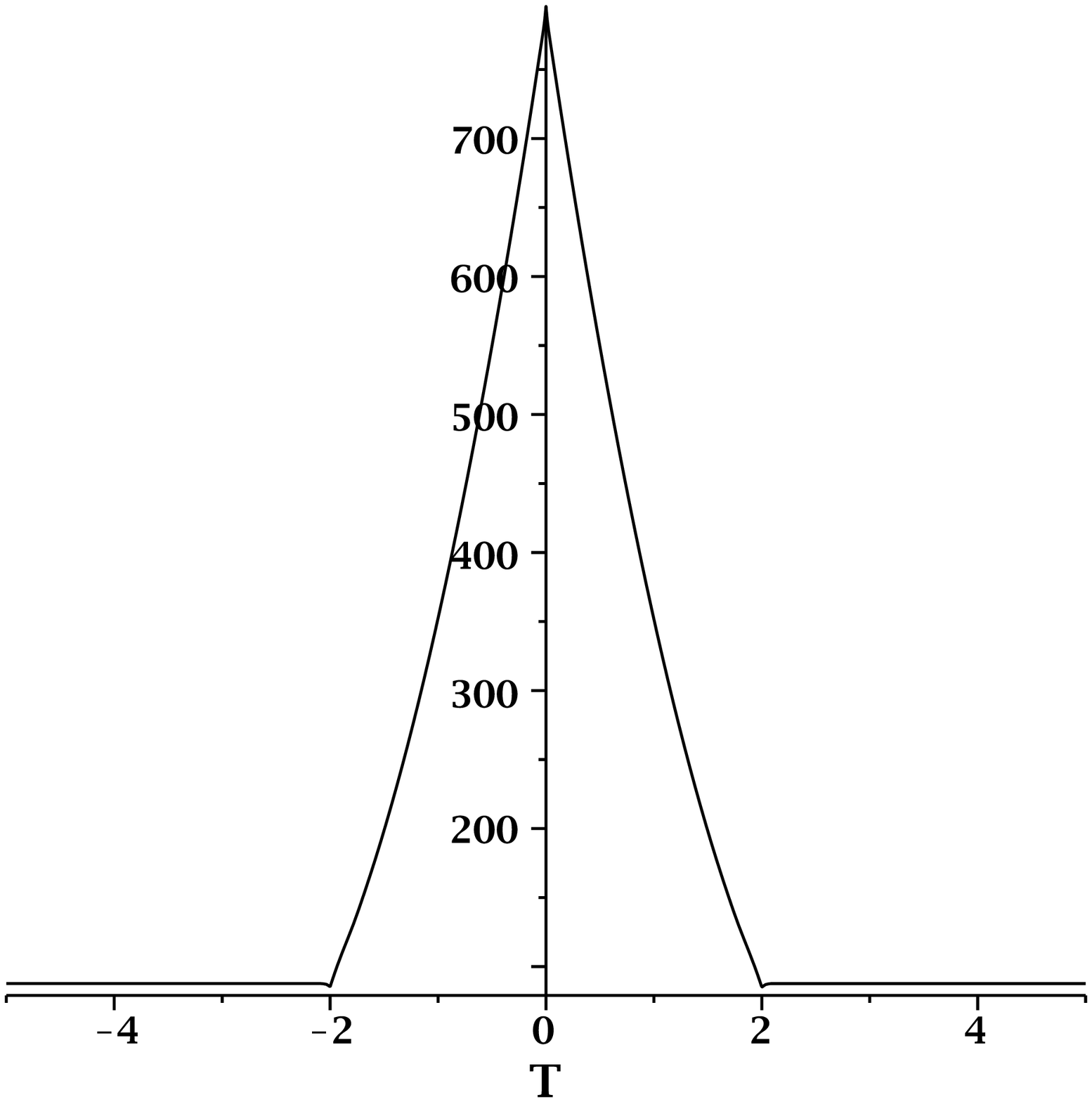}}}
\caption{Verification of the Voigt-type inequality for the effective properties of the nonlinear composite.}%
\label{minor2}
\end{center}
\end{figure}

The respective results are presented in Fig. \ref{minor1},\ref{minor2}.

\begin{figure}[h!]
\begin{center}
\resizebox*{14cm}{!}{{\includegraphics{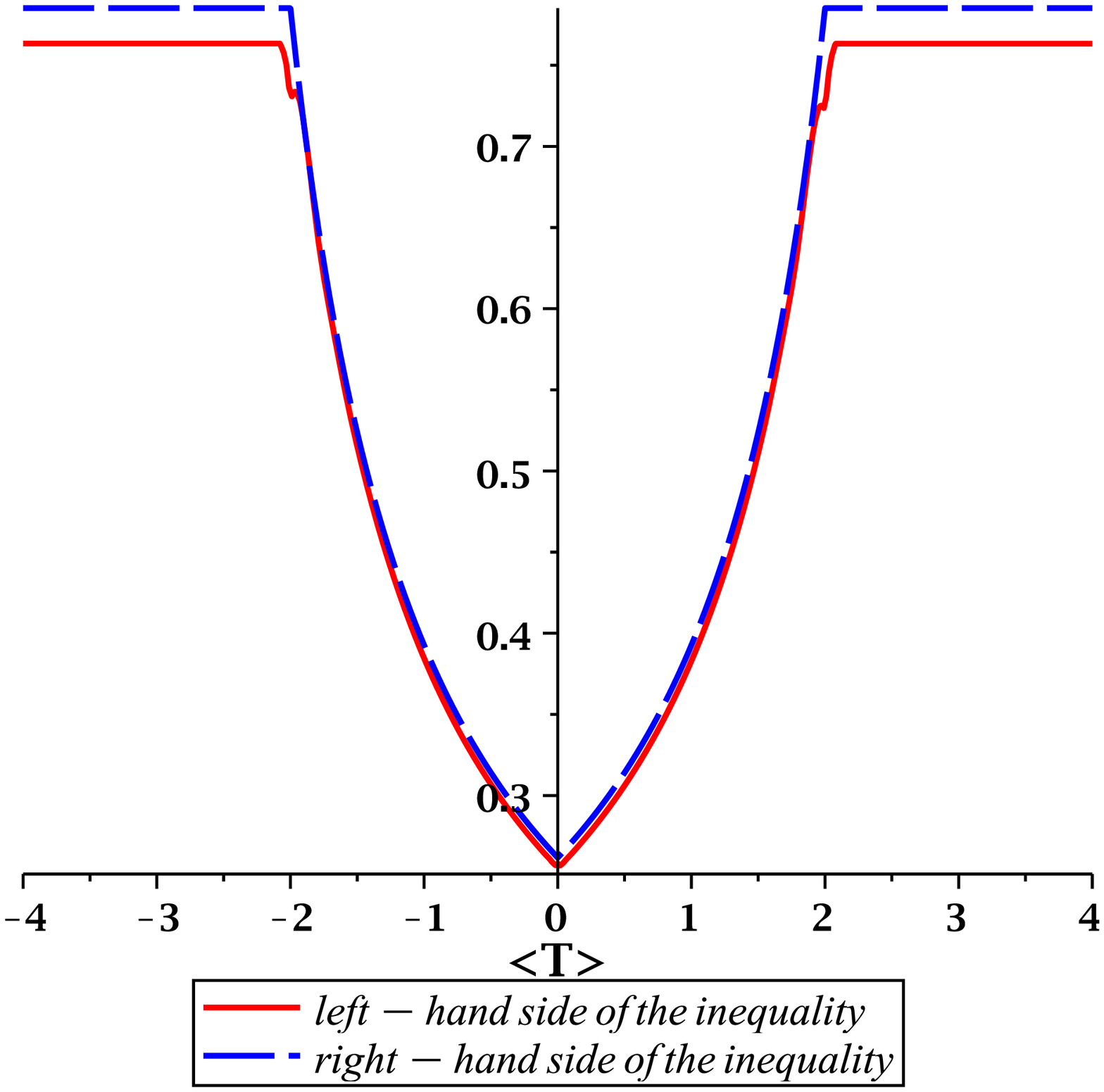},\includegraphics{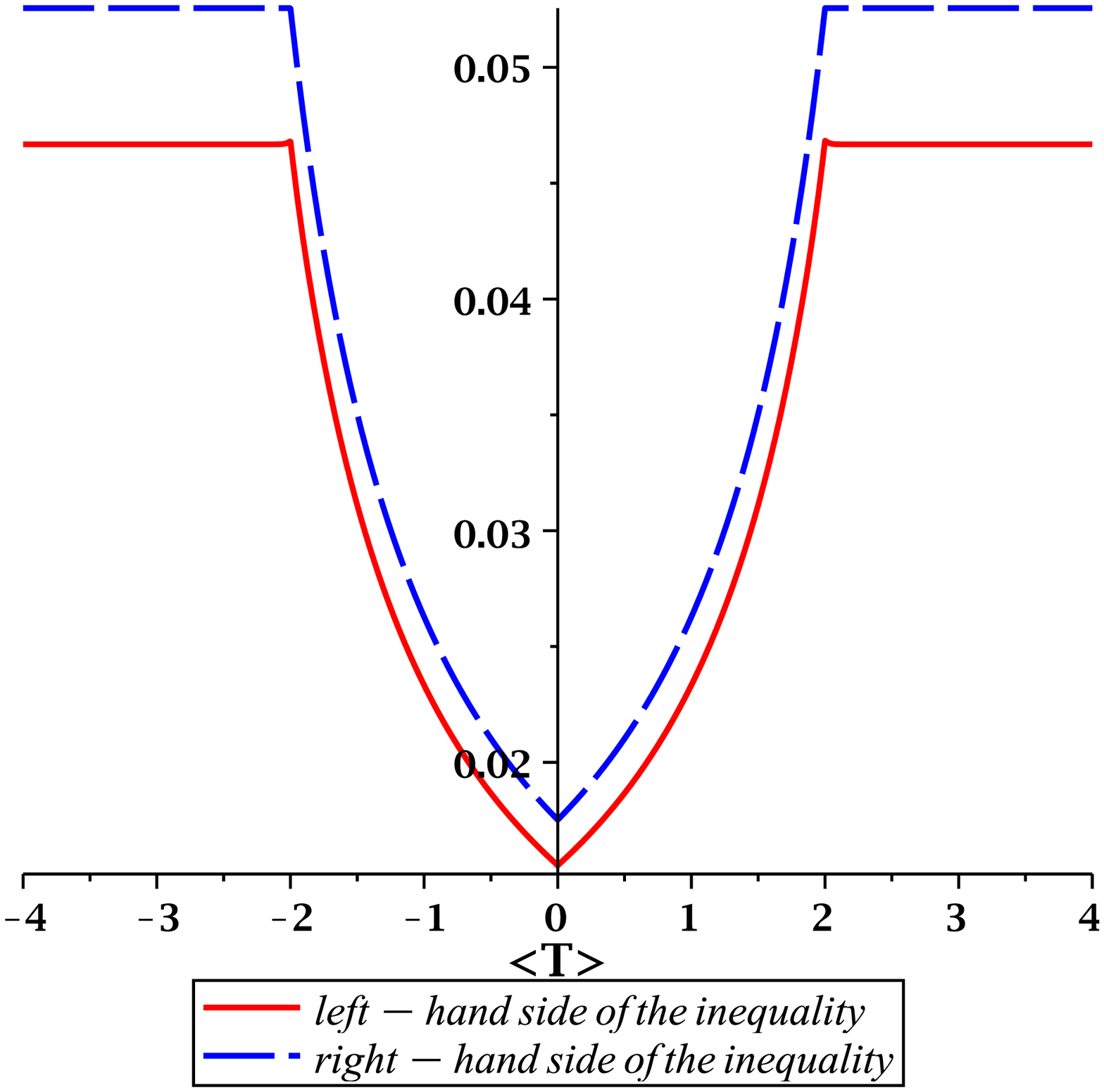}}}
\caption{Verification of the Hashin-Shtrikman bounds (\ref{HS1}) -- Fig \ref{SH} a) and (\ref{HS2}) -- Fig \ref{SH} b).
Solid (dash) lines correspond to the left (right)-hand sides of the inequalities.
}%
\label{SH}
\end{center}
\end{figure}

Now we check feasibility of the Hashin-Shtrikman bounds extended in \cite{GalTelTok} for the case of quasi-linear composite. These estimates
are more narrow than the elementary bounds (\ref{el_bounds}) and can be written in our notation (comparing with \cite{GalTelTok}):
\begin{equation}\label{HS1}
tr \left[ (\Lambda_e(T)-\lambda(T)I)^{-1}\right]\leq \frac{1}{\mu_2(T)-\lambda(T)}+\frac{1}{\mu_1(T)-\lambda(T)},
\end{equation}
and
\begin{equation}\label{HS2}
tr \left[ (\lambda_k(T)I-\Lambda_e(T))^{-1}\right]\leq \frac{1}{\lambda_k(T)-\mu_2(T)}+\frac{1}{\lambda_k(T)-\mu_1(T)},
\end{equation}
where $tr A=A_{jj}$, ($j=1,2$).
The left and right hand sides of the inequalities  (\ref{HS1}), (\ref{HS2}) are presented on Fig. \ref{SH},
where solid (dash) lines correspond to the left (right) -hand sides of the inequalities.

\newpage

\subsection{Comparison with the results for random composite and discussions}

According to \cite{MishSev}, the thermal effective conductivity $\lambda_e$ of an isotropic composite with temperature dependent and proportional conductivities
of the components may be computed by the standard homogenization techniques in the following form:
\begin{equation}\label{MS}
\Lambda_e(T)= \lambda(T)\cdot \Lambda_e ,
\end{equation}
where $\Lambda_e$ is the tensor of effective conductivity of a linear problem with the same ratio $C_k$ between the conductivity of the matrix and the inclusions as we have for the nonlinear one (see (\ref{main_cond})). Thus, in order to find the effective conductivity for a such kind of composites it is sufficient to find only the effective conductivity of corresponding linear problem. In this particular case we found:
\begin{equation}
\label{Lambda1}
\Lambda_e=\left(\begin{array}{cc}
1.524131 & 0.000027\\
0.000027 & 1.650632
\end{array}\right),
\end{equation}
with the same accuracy $10^{-6}$ as discussed above.

\begin{figure}[h!]
\begin{center}
\resizebox*{14cm}{!}{{\includegraphics{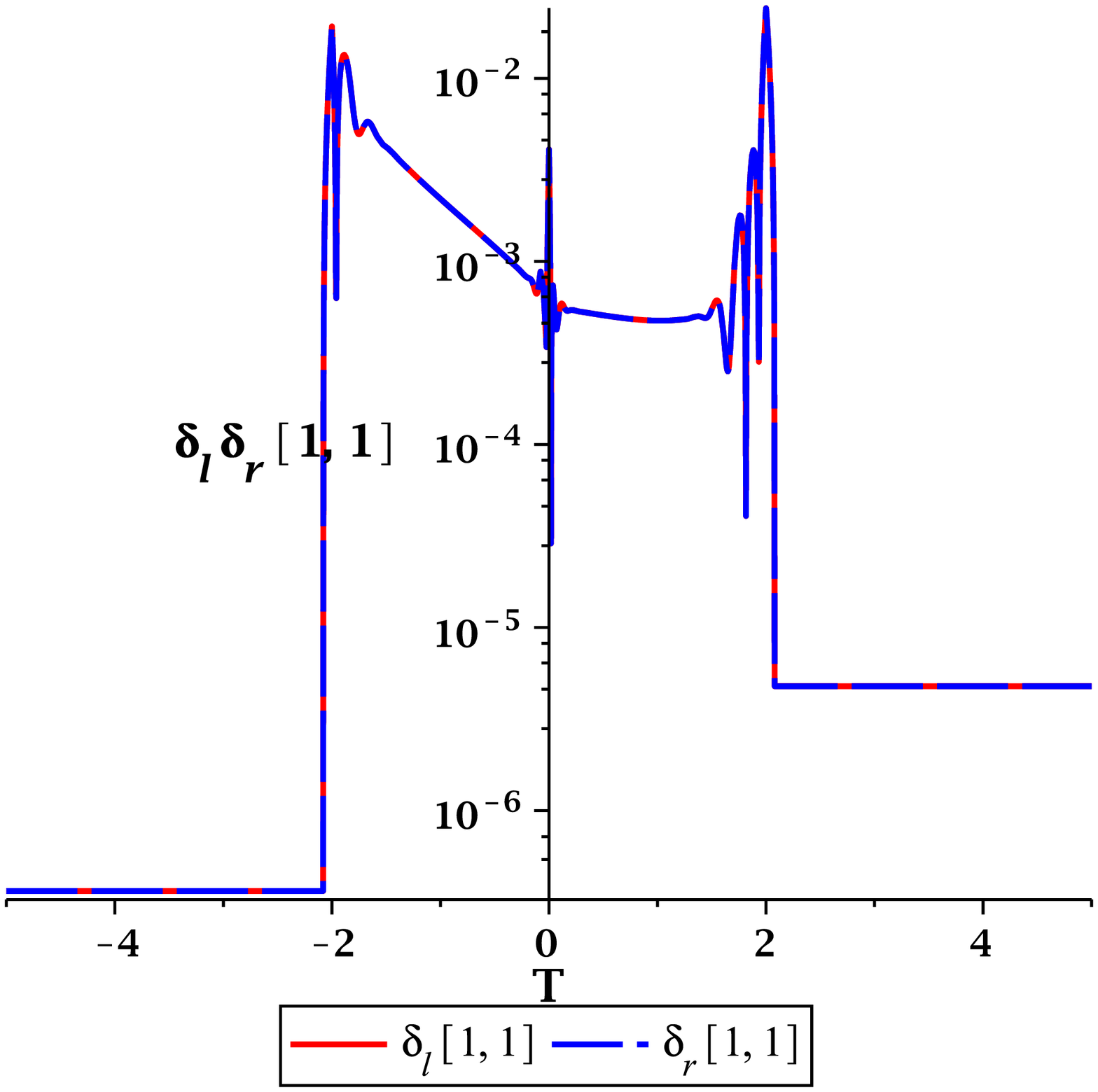},\includegraphics{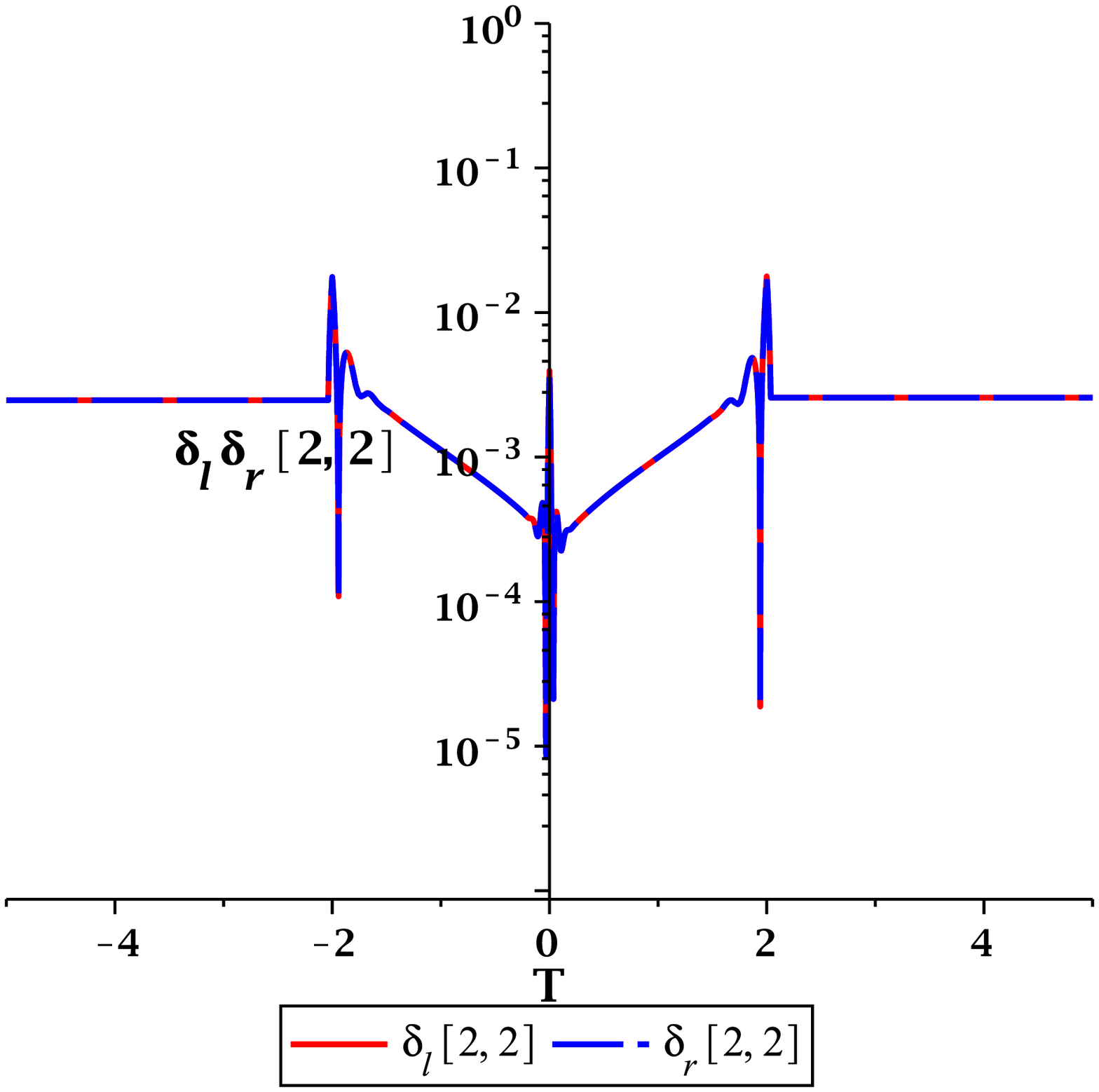}}}
\caption{The main diagonal components ([1,1] and [2,2]) of the tensors $\delta_l$ and $\delta_r$.}%
\label{dif11_log}
\end{center}
\end{figure}

\begin{figure}[h!]
\begin{center}
\resizebox*{14cm}{!}{{\includegraphics{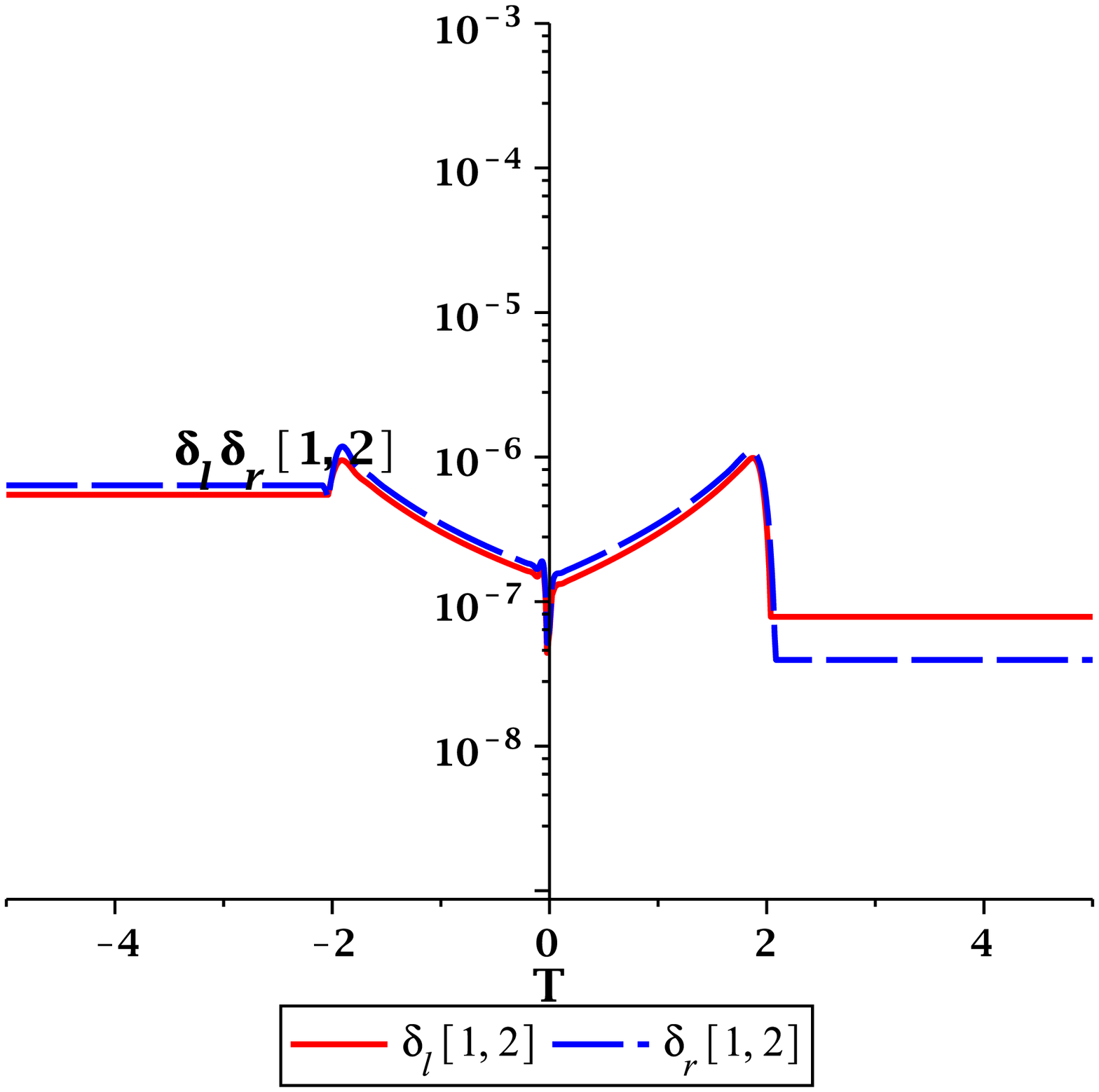},\includegraphics{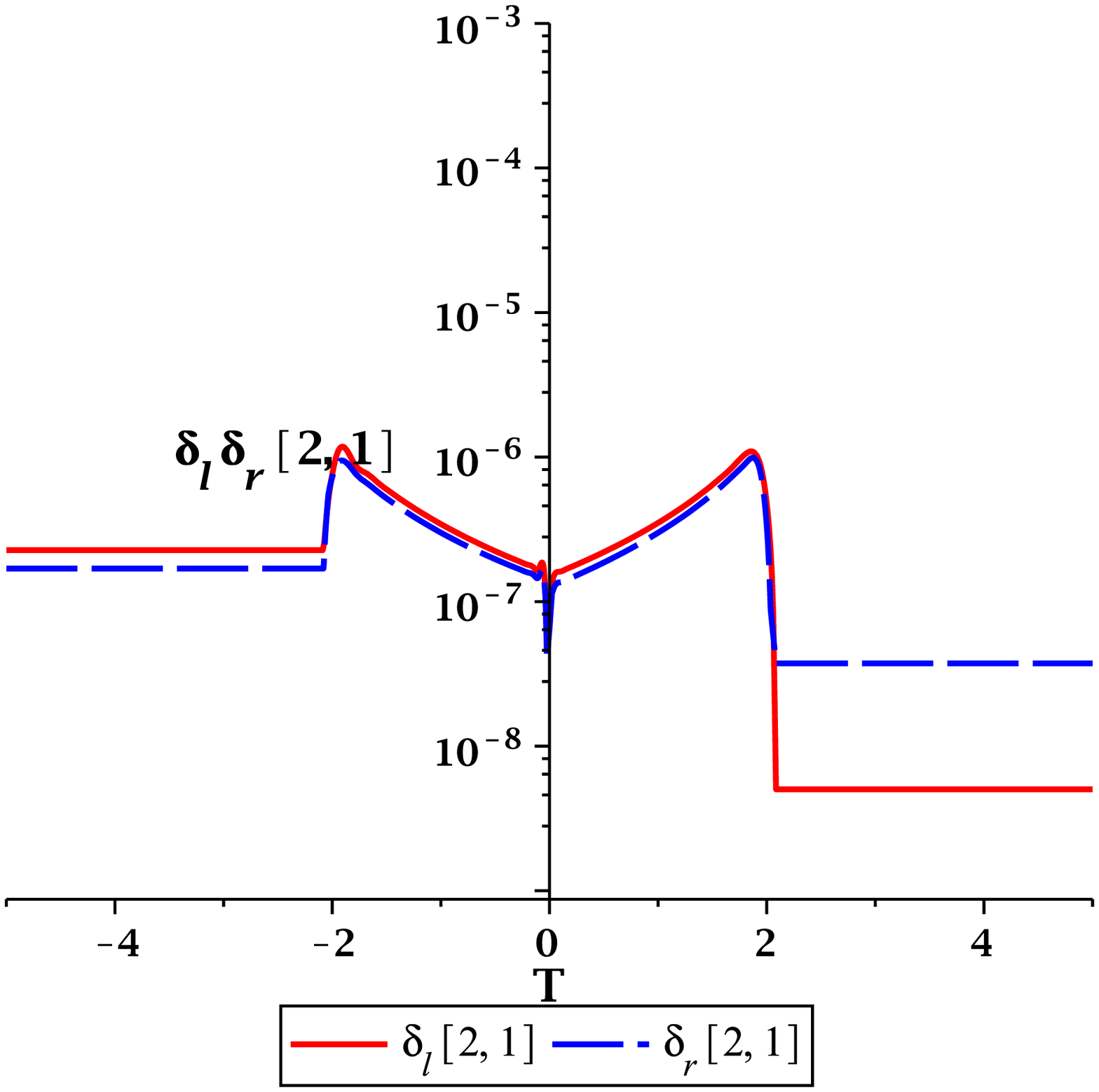}}}
\caption{Deviation in the components [1,2] and [2,1] of the tensors $\delta_l$ and $\delta_r$. }%
\label{dif22_log}
\end{center}
\end{figure}

We check whether the relationship (\ref{MS}) is useful for the doubly periodic composite considered in this paper which is not isotropic
(see Fig. \ref{fig:4_inclusions}).
For this, we calculate two relative errors
(bearing in mind that the computed values are tensors):
\begin{equation}
\label{tenz_komp_comp}
\begin{array}{c}
\delta_l=(\Lambda_e (T)-\lambda(T)\cdot\Lambda_e)\cdot(\Lambda_e (T))^{-1}, \\[3mm]
\delta_r= (\Lambda_e (T))^{-1}\cdot(\Lambda_e(T)-\lambda(T)\cdot\Lambda_e).
\end{array}
\end{equation}
The components of the tensors $\delta_l$ and $\delta_r$ from (\ref{tenz_komp_comp}) are represented on Fig. \ref{dif11_log},\ref{dif22_log}, respectively.
The curves are given in the logarithmic scale to clearly indicate the difference. By the dash line we present the components of $\delta_r$ while for the solid line corresponds to the components of $\delta_l$. Although the model used in \cite{MishSev} is less accurate than the model analyzed in this paper, our computations show perfect correlation between the results. The largest deviations (near $2\%$) take please
near the points whether the original functions are not smooth. This difference is observed only for the components situated in the tensor on the main diagonal.
The other two components are rather identical taking into account the computational accuracy. The last result is the direct consequence of the law anisotropy of the problem under consideration.


 \vspace{3mm} {\bf Acknowledgements.} G. Mishuris is grateful for support from the FP7
IRSES Marie Curie grant TAMER No 610547.
D. Kapanadze is supported by Shota Rustaveli National Science Foundation with the grant number 31/39. E. Pesetskaya is supported by Shota Rustaveli National Science Foundation with the grant number 24/03.

\end{document}